\documentclass[fleqn,12pt]{article}
\usepackage{graphicx}
\usepackage{epsf}
\usepackage{amssymb}
\usepackage{amsmath}
\usepackage{fancyheadings}

\newcommand{\be}{\begin{equation}}
\newcommand{\ee}{\end{equation}}
\newcommand{\bea}{\begin{eqnarray}}
\newcommand{\eea}{\end{eqnarray}}
\newcommand{\bean}{\begin{eqnarray*}}
\newcommand{\eean}{\end{eqnarray*}}

\newcommand{\Cov}[2]{\mbox{Cov}\BRK{#1,#2}}

\newcommand{\bfs}{{\mathbf s}}

\newcommand{\bfb}{{\mathbf b}}
\newcommand{\bfu}{{\mathbf u}}
\newcommand{\bfz}{{\mathbf z}}

\newcommand{\brk}[1]{\left(#1\right)}          
\newcommand{\Brk}[1]{\left[#1\right]}          
\newcommand{\BRK}[1]{\left\{#1\right\}}        
\newcommand{\Average}[1]{\left<#1\right>}      
\newcommand{\Abs}[1]{\left| #1 \right|}        

\newcommand{\pd}[2]{\frac{\partial#1}{\partial#2}}
\newcommand{\deriv}[2]{\frac{d#1}{d#2}}

\newcommand{\half}{\frac{1}{2}}

\newtheorem{lemma}{Lemma}[section]

\newtheorem{assumption}{Assumption}[section]

\begin{document}

\title {Optimal Prediction for Hamiltonian partial differential equations}
\author{Alexandre J. Chorin\footnotemark[2] 
\and  Raz Kupferman\footnotemark[3] 
\and Doron Levy\footnotemark[4]}
\date{}

\renewcommand{\thefootnote}{\fnsymbol{footnote}}
\footnotetext[2]
{Department of Mathematics, University of California, 
Berkeley, CA 94720, and Lawrence
Berkeley National Laboratory; {\tt chorin@math.berkeley.edu}}
\footnotetext[3]
{Institute of Mathematics, The Hebrew University, 
Jerusalem 91904, Israel; {\tt razk@math.huji.ac.il}}
\footnotetext[4]
{Department of Mathematics, University of California, 
Berkeley, CA 94720, and Lawrence
Berkeley National Laboratory; {\tt dlevy@math.berkeley.edu}}
\renewcommand{\thefootnote}{\arabic{footnote}}
\maketitle

\begin{abstract}
Optimal prediction methods compensate for a lack of resolution in the
numerical solution of time-dependent differential equations through
the use of prior statistical information.  We present a new derivation
of the basic methodology, show that field-theoretical perturbation
theory provides a useful device for dealing with quasi-linear
problems, and provide a nonlinear example that illuminates the
difference between a pseudo-spectral method and an optimal prediction
method with Fourier kernels.  Along the way, we explain the
differences and similarities between optimal prediction, the
representer method in data assimilation, and duality methods for
finding weak solutions. We also discuss the conditions under which a
simple implementation of the optimal prediction method can be expected
to perform well.

\end{abstract}



\bigskip
\noindent
{\bf AMS(MOS) subject classification:} Primary 65M99; 
secondary 81T15, 35Q55, 65M70.



\section{Introduction}			\label{sec:introduction}
\setcounter{equation}{0}
\setcounter{figure}{0}
\setcounter{table}{0}

We consider a Cauchy problem
\begin{equation}
\left\{ 
\begin{array}{l} 
u_t = R\left(x,u, u_{x}, u_{xx},\ldots \right),
\label{eq:basic} \\ \\ 
u(x,0) = u_{0}(x), 
\end{array}
\right.
\end{equation}
where $R(x,u,u_x...)$ is a (generally nonlinear) function of
$u=u(x,t)$, of its spatial derivatives, and of the independent
variable $x$, in any number of dimensions; subscripts denote
differentiation.  We assume that we cannot afford to use enough
computational elements (for example, mesh points) to resolve the
problem adequately, or that we do not wish to use many computational
elements because we are only interested in some of the features of the
solutions and it seems to be wasteful to compute all the details.  We
also assume that we have some prior statistical information about the
distribution of possible solutions.  The question we address is how
the prior statistics can be used to obviate the need for resolution.
In contrast to problems that arise in some applications, especially in
geophysics \cite{bennett}, we assume that the equations we are solving
are fully known, and that the data are knowable in principle, even if
we may not be able or willing to store them all in a computer's
memory.

In the present paper we assume that our statistical information
consists of an invariant measure $\mu$ on the space of solutions,
i.e., we assume that the initial data are sampled from a probability
distribution on the space of data, and that, in principle, if one
takes one instance of initial data after the other and computes the
solutions produced by (\ref{eq:basic}) at any later time $t$, then the
set of solutions obtained at that later time $t$ (viewed as functions
of the spatial variables) has the same probability distribution as the
set of initial data. We assume further that this probability
distribution is explicitly known. In section 4 below we give an
example where these assumptions are satisfied and provide more precise
definitions.  Our present assumptions may be unnecessarily strong for
some practical problems, but they simplify the exposition and it is
often easy to weaken them.  We shall call the invariant measure the
{\em prior measure}; thus the prior measure is the distribution of the
data before anything has been specified about a particular problem.
In Hamiltonian systems a natural prior measure is the canonical
measure induced by the Hamiltonian ${\cal H}$, i.e., a measure defined
by the probability density
\[
f(u) = Z^{-1} e^{-{\cal H}(u)/T},
\]
where $Z$ is a normalization constant and the parameter $T$, which determines the variance of the density, is known as the 
``temperature''.  We do not assume that the differential equation
(\ref{eq:basic}) admits a unique invariant measure; in cases of
non-uniqueness the right choice of measure is part of the formulation
of the problem.

We further assume that all we know, or care to know, at time $t=0$, is
a small set of data.
In the present paper we choose these data to be of the form
\begin{equation}
(g_\alpha,u) = \int g_{\alpha}(x) u(x,0) \,dx = V_{\alpha}, 
\qquad \alpha=1,\dots,N,   
\label{eq:eq2}
\end{equation}
where the $g_{\alpha}=g_{\alpha}(x)$ are suitably chosen kernels (the
choice of kernels $g_{\alpha}$ is at our disposal).  The question we
are asking can now be rephrased as follows: How do we best predict the
future using only $N$ variables defined as in (\ref{eq:eq2}), given
the prior measure $\mu$?

At time $t=0$ one can, at least in principle, find the mean solution
conditioned by the data (\ref{eq:eq2}), for every point $x$, i.e., at
each $x$ one can average over all those functions in the support of
the invariant measure that also satisfy the conditions (\ref{eq:eq2})
and find the mean $v$ of $u$ given the values $V_{\alpha}$,
$\alpha=1,...,N$; symbolically,
\begin{equation}
v(x) = E\Brk{u(x)| V_1,...,V_N},  
\label{eq:const.a}
\end{equation}
where $E[\cdot]$ denotes an expectation value. This is a regression
problem, which can be solved by standard tools (\cite{papoulis}, see
also below).  If one is given only the information contained in
(\ref{eq:eq2}), this regression is the right substitute for a more
detailed knowledge of the data.  To perform the regression at a later
time $t$ we need appropriate conditions at that later time, which
should encode the later effect of the initial data (\ref{eq:eq2}).
The problem now at hand is how to find these conditions and how to do
it efficiently, in linear and in nonlinear problems.

We have already explained our basic approach in
\cite{ckk1,ckk2,ckk3}.  In the present paper we explain
it in a different way which we hope is more transparent; the
suggestion that a weak formulation is the right starting point is due
to Gottlieb \cite{gottlieb}.  The paper has several goals: to explain
what constitutes the novelty of our approach, to show how it meshes
with field-theoretical perturbation methods, to give a simple,
explicit, nonlinear example of the ways in which the basic algorithm
differs from, and is superior to, a numerical method that uses the
information in the initial data without taking advantage of the prior
statistics, and to discuss the domain of applicability of the current
implementation of the basic ideas.  Our examples are of Schr\"odinger
type.

Other work along somewhat analogous lines includes Scotti and Meneveau
\cite{scotti}, where a ``fractal''interpolation can be viewed as an
analog of sampling a measure, and a stochastic construction by
Vaillant \cite{vaillant}.

A comment on notations: The notation in equation (\ref{eq:const.a}) is
clear but cumbersome.  A shorter version is: $v=E[u|V]$. We shall also
use the simpler but less transparent physicists' notations:
$\left<u\right>$ for $E[u]$ and $\left<u\right>_V$ for $E[u|V]$.  Both
notations $\left<u\right>_V$, $E[u|V]$ are of course shorthands for
$\left<u\right>_{V_{\alpha},g_{\alpha}}$ etc., since the conditioning
(\ref{eq:eq2}) depends on the kernels as well as on the right-hand
sides.

The paper is organized as follows: In section 2 we present an overview
of optimal prediction. In sections 3, 4, we present some needed
background material on regression and invariant measures. In sections
5, 6 we explain how perturbation theory can be used to implement
optimal prediction in nonlinear problems, and provide examples of
success and failure in problems with sparse data.  In section 7 we
provide a detailed analysis of a problem in which perturbation theory
is carried only to zero-th order (but with a ground state that our
previous work shows to be optimal), and in which the kernels are
trigonometric functions.  These simplifications yield results that are
particularly transparent.  Conclusions are drawn in a final section.


\section{Weak solutions, regression, and prediction}         
\label{sec:weak}

\setcounter{equation}{0}
\setcounter{figure}{0}
\setcounter{table}{0}


We start with well-known considerations about weak solutions of linear
equations; they lead naturally to an algorithm for solving
underresolved linear equations, which will be useful below, in
particular because we shall be able to present our proposals by way of
a contrast. Consider the equation:

\begin{equation}
u_t = L u,  
\label{eq:linear}
\end{equation}
where $L$ is a linear operator. Multiply (\ref{eq:linear}) by a smooth
test function $g$; for simplicity assume that the boundary conditions
on $u$ are periodic, and thus $g$ can be also assumed
periodic. Integration  over a periodic domain in $x$
and between $t=0$ and $t=\tau$, followed by an integration by parts,
results in
\begin{equation}
\int_{t=0}^{\tau}\int_{x} \brk{g_t + L^{\dagger} g} u \,dx \,dt + 
(u(\tau),g(\tau)) - (u(0),g(0)) = 0, 
\label{eq:weak}
\end{equation}
where $(u(t),g(t))$ denotes $\int u(x,t)g(x,t) \,dx$ and $L^{\dagger}$
is the adjoint of $L$.  A weak solution of (\ref{eq:linear}) is a
function $u(x,t)$ that satisfies (\ref{eq:weak}) for all test
functions $g$ (see e.g.\cite{lax}).  In particular, if the test
functions $g_{\alpha}$ satisfy the adjoint equation
\begin{equation}
\frac{\partial g_{\alpha}}{\partial t} + L^{\dagger} g_{\alpha}=0,
\qquad \alpha=1,\dots,N,
\label{eq:adjoint}
\end{equation}
then $u$ is a weak solution of (\ref{eq:linear}) provided
$(u,g_{\alpha})$ is a constant independent of $t$ for each $\alpha$.
This observation produces a possible numerical method for solving
(\ref{eq:linear}): One can construct a collection of functions
$g_{\alpha}$ that satisfy the adjoint equation (\ref{eq:adjoint}),
find the numerical values $V_{\alpha}$ of the inner products
$(u,g_{\alpha}), (g=g(x,0))$, $u$ being the initial data, and finally
reconstruct the weak solution at a later time $\tau$ from its inner
products with the functions $g(x,\tau)$.  These functions, $g(x,t)$,
can be found at time $t=\tau$ if they are known at $t=0$.

If this method is used with a small number of functions (i.e., small
$N$), the solution at a time $t>0$ will be underdetermined. In this
case one can use the invariant measure $\mu$ to ``fill in'' the gaps
through regression, i.e., replace the weak solution which is not
completely known by its average (as determined by $\mu$) over all
solutions that satisfy the $N$ conditions $(g_{\alpha},u)=V_{\alpha}$.
In other words, replace the function $u$ which is not completely known
by the regression $v=v(x,t)=\left<u\right>_V$.  Some technical
background on regression follows in \S\ref{sec:regression}.

Note that so far, the construction resembles what is quite commonly
done in underdetermined linear problems (for example, in the context
of data assimilation \cite{bennett}; the function $g_{\alpha} $ are
analogous to the ``representers'' which are used there).  The
construction is computationally useful in certain linear problems even
when alternate ways of finding future regressions are available, as
happens whenever time evolution and averaging commute (see
\cite{ckk2}). However, the construction just presented is restricted
to linear problems, and the amount of work needed to evaluate and
store the kernel functions $g_{\alpha}$ may not be trivial.  It is
also clear that not all of the available information has been used, as
the evolution of the kernels $g_{\alpha}$ is independent of the
invariant measure $\mu$, and the measure used in the regression need
not be connected with the differential equation.  Indeed, in
geophysical applications the measure is chosen according to
considerations quite extraneous to the differential equations which
may be only partly known.

We now wish to show how the machinery can be modified through the use 
of the invariant measure in the evolution equations for the 
conditions, so that it becomes more efficient as well as generalizable 
to nonlinear problems.

First, we must view the calculations differently. The measure $\mu$
induces a measure on the space of the data $u$ that satisfy the
initial conditions: the {\em conditional measure\/} $\mu_V$.
\emph{The conditional measure is not invariant}; for example, if the initial
conditions consist of $N$ point-values of the functions $u(x,0)$,
i.e., if we assume that at $t=0$ the functions satisfy
$u(x_{\alpha},t=0)=V_{\alpha}$ (the functions $g_{\alpha}$ are then
$\delta$ functions), there is no reason to believe that these
conditions will be satisfied at all subsequent times with the same
$V_{\alpha}$ by the solutions of the differential equations that arise
out of these initial data.

Imagine first that one can sample the initial conditional measure
$\mu_V$; find the (in general, weak) solution of equation
(\ref{eq:linear}) that has this datum as initial condition, and
perform this procedure repeatedly. At time $t$ this produces an
ensemble of functions $u(x,t)$ which inherits a measure from the
initial data.  In principle, this measure is well-determined; we wish
to determine it in practice, and then to average with respect to this
measure so as to obtain what we call an {\em optimal prediction}.
Note that if the temporal evolution governed by (\ref{eq:basic}) is
ergodic with respect to the invariant measure, the conditional measure
will eventually relax to the invariant measure and the initial
conditions will be forgotten; we are in fact dealing with a
computational analog of non-equilibrium statistical mechanics.

Given the initial conditional measure, we can find the
statistics of $Lu$, or, in the general case (\ref{eq:basic}), the
statistics of $R(u)$ and consequently the statistics of $u_t$; thus,
the evolution of the measure $\mu_V$ can be determined for a
short time interval $\Delta t$.  We cannot go beyond a short time
interval because the measure $\mu_V$ at time $\Delta t$ can no longer
be described as the invariant measure conditioned by the conditions
(\ref{eq:eq2}), at least not with the same functions $g_{\alpha}$ and
the same $V_{\alpha}$.  This leads us to the closure assumption:

\begin{assumption}[Closure]
The conditional probability measure at time $t$, $\mu_V(t)$, can be
approximated by 
\begin{equation}
\mu_V(t) = \mu_{V(t)},   
\label{eq:closure}
\end{equation}
where the left hand side $\mu_V(t)$ is the measure conditioned by
the initial data (\ref{eq:eq2}), while the right hand side is the
invariant measure conditioned by $N$ affine conditions of the form
\[ 
(g_{\alpha}(t),u(t)) =
\int g_{\alpha}(x,t) u(x,t) \,dx = V_{\alpha}(t).
\] 
\end{assumption}
The kernels $g_{\alpha}(t)$ and the values $V_{\alpha}(t)$ of the
inner products will generally be different at time $t$ than at time
$t=0$, but it is assumed here that the affine form of the conditions
and their number remain constant.  We have already seen that in the
linear case equation (\ref{eq:closure}) is in fact a theorem (a
different analysis of this fact was given in \cite{ckk2}).

It should be noted that even though this is the assumption that will be
used in the present paper, it may be unduly restrictive in other
situations; there is nothing magical about keeping the number of
conditions fixed, and the conditions need not in general be affine.

What we wish to do is to advance the measure to time $\Delta t$ and
then to find the conditions that produce that new measure at the new
time, when these conditions condition the invariant measure.  In that
way, the whole process can be taken a step further and then repeated
as often as one may wish.  In the linear case we have already produced
a recipe for updating the conditions: there, we let $g_{\alpha}$
satisfy the adjoint equation and keep the numbers $V_{\alpha}$ fixed.
We now propose different, {\em approximate}, ways of finding
conditions that describe the evolving measure.  We are going to do so
by matching moments; on one hand, we will calculate moments of the
conditional measure by regression from the old moments, and on the
other hand, we will produce conditions that produce the new moments
from the invariant measure; this will produce equations of motion for
the conditions.

More specifically, suppose that we compute 
$Nq$ moments of $u$, at time $\Delta t$.
For example, set $q=2$, and compute the means and the variances
with respect to the conditioned measure of the random variables
$u(x,\Delta t)$ at each of the $N$ points $x_{\alpha}$.  We can do
this knowing the conditions at time $t=0$ and the invariant measure.
On the other hand, suppose we let the functions $g_{\alpha}$ depend on
$q-1$ parameters; for example, if $q=2$, then we can pick
$g_{\alpha}(x)=\exp((x-x_\alpha)^2/\sigma_{\alpha})$, where the
``centers'' $x_{\alpha}$ are fixed and the numbers $\sigma_{\alpha}$
may be allowed to vary in time and will serve as our parameters.  If
we write down the requirement that the moments we calculated at
time $\Delta t$ match the moments produced by conditioning the
invariant measure by affine conditions  with unknown values of
parameters such as the $\sigma_{\alpha}$ and unknown values of the
right-hand sides $V_{\alpha}$, we obtain $Nq$ algebraic equations for
the parameters, which we can try to solve.  If we solve these equations
successfully, we obtain a set of simultaneous \emph{ordinary}
differential equations for the parameters and the moments.

Before carrying out such a calculation, two remarks: There is no a
priori guarantee that the algebraic equations we will obtain can be
always solved: The basic assumption may fail, and the choice of
parameters may be unsuitable. A good numerical program will inform us
that a solution cannot be found.  However, if a solution is found, the
resulting moments are \emph{realizable}.  It is well-known that
closures may well produce moments that not only fail to solve the
problem at hand but do not solve any problem, because there is no
stochastic process that admits the computed moments as its moments
(see e.g.\cite{proudman}).  This is often a major difficulty in the
formulation of mean equations, and it does not arise here.

We limit ourselves here to the simplest case with
$q=1$, i.e., for each $\alpha, \alpha=1,...,N$ we keep track of a
single quantity, which we choose to be the mean value
of the inner product of $g_{\alpha}$ and $u$, 
$\left<(g_{\alpha},u)\right>_V= (g_{\alpha},\left<u\right>_V)$, 
while at the same time we
modify a single parameter in each condition; that single parameter is
chosen to be the value $V_{\alpha}$ of the $\alpha$-th product.  Thus
we must have for $\alpha=1,...,N$,
\[
\frac{d}{dt}V_{\alpha} =
\Average{(g_{\alpha},u_t)}_V =
\Average{(g_{\alpha},R(u))}_V =
(g_{\alpha},\Average{R(u)}_V),
\]
where $R(u)$ is the right-hand side of equation (\ref{eq:basic}).  
The final equation,
\begin{equation}
\frac{dV_{\alpha}}{dt} = (g_{\alpha},\Average{R(u)}_V),   
\qquad \alpha = 1 \ldots N,
\label{eq:final.main}
\end{equation}
is the main equation used in the present paper. The more transparent
mathematical notation,
\[
\frac{dV_{\alpha}}{dt} = E\Brk{(g_{\alpha},R(u))|V_1,\dots,V_N},
\qquad \alpha = 1 \ldots N,
\]
makes explicit the fact that we are dealing with a system of $N$
ordinary differential equations. Once the $V_{\alpha}$ are found at
time $t$, a regression can be used to find the average solution at any
point $x$ (see \S\ref{sec:regression}).  We shall call an algorithm
that uses, in the matching of moments and parameters, only means of
the unknown solutions and no higher moments, a first-order prediction
scheme; in the present paper this is the only prediction scheme we
shall use.  We hope to demonstrate that a first-order prediction is
often an improvement over algorithms that take no cognizance of the
invariant measure, in the sense that it requires less computational
labor than the alternatives; higher-order and more sophisticated
optimal prediction schemes will be described in subsequent work.

It is useful to contrast our algorithm with the one at the beginning
of the section: We are keeping the kernels $g_{\alpha}$ fixed 
while changing the values $V_{\alpha}$ of the conditions, while the
``natural'', linear construction at the beginning of the section did
the opposite. 

When can we expect the first order optimal prediction scheme to be
accurate? In a linear problem $u_t=Lu$, if the kernels $g_{\alpha}$
are eigenfunctions of the operator $L^\dagger$ adjoint to $L$, then
one can readily see from the analysis at the beginning of the section
that equations (\ref{eq:final.main}) are exact, and this remains true
if the $g_{\alpha}$ span an invariant subspace of $u_t+L^{\dagger}=0$.
Lowest-order optimal prediction should be accurate as long as the
$g_{\alpha}$ span a space that is approximately invariant under the
flow induced by $L^\dagger$ (see
\cite{gottlieb},\cite{hald:optimal}).  This remark also provides a
recipe for choosing the kernels.  Something similar remains true in
nonlinear problems \cite{evans}, \cite{mori},\cite{zwanzig}: Define
the space of functions spanned by the kernels $g_{\alpha}$ to be the
resolved part of the solution; equations (\ref{eq:final.main}) should
yield an accurate prediction of this resolved part (including a
correct accounting for the effect of unresolved components on the
resolved components) as long as there is no substantial transfer of
information from the resolved part to the unresolved part and back; if such
information 
transfer should occur, a correct description of the flow should
include an additional, ``memory'', non-Markovian term.  This remark
also points out that equations (\ref{eq:final.main}) should not be
accurate for very long times (because memory terms are important to
the description of decay to equilibrium), and should be better at low
temperatures $T$ than at high temperatures (because the decay to
equilibrium should be more rapid at high $T$).

There are other ways to ensure the accuracy of a first order optimal
prediction scheme: The quantities $(g_{\alpha},R(u))$ are, of course,
random variables whose distribution depends on the measure
$\mu_V$. If the standard deviations of these variables are small,
then equations (\ref{eq:final.main}) are be good approximations to
the exact solution.  Indeed, equation (\ref{eq:final.main}) merely
equates these random variables to their means; a higher order
approximation would take into account the variance of these variables
as well. The smaller the variance of the variables
$(g_{\alpha},R(u))$, the smaller the error we expect; to some extent,
we can control this variance by choosing the kernels
appropriately. The larger the support of the kernels, the more these
variables represent spatial averages, and the slower we may expect
their variance to grow; thus if the scheme is to be accurate, it is
most likely that the functions $g_{\alpha}$ should not be narrow,
$\delta$-function-like objects.  Furthermore, the kernels should not
have disjoint supports (for an analogous observation in the theory of
vortex methods, see, e.g., \cite{chorin2}).  Rigorous error bounds for
the linear case can be found in \cite{hald:optimal}.  These conditions
allow us to call the variables $(g_\alpha,u)$ that appear in
(\ref{eq:eq2}) {\em collective variables}; they are groupings of
variables.  Note that as the number $N$ of collective variables
increases, optimal prediction homes in an ever smaller set of initial
data, and the variance of the variables $(g_{\alpha},R(u))$ should
decrease.

Finally, the presentation above started from a discussion of weak
solutions, and indeed in all the examples below the solutions will be
weak; why is that so? We shall present a detailed mathematical
analysis elsewhere; here it should suffice to comment that if the
solutions are not highly oscillatory on several levels, there is less
interest in analyzing methods that fail to resolve them; solutions
that do oscillate significantly on several scales appear, on the
largest scale, as non-smooth and therefore weak.


\section{Regression for Gaussian variables}    \label{sec:regression}

\setcounter{equation}{0}
\setcounter{figure}{0}
\setcounter{table}{0}

It was mentioned in the previous section that once the data that
condition the state of the system are found, or, at time $t=0$ once
the initial data have been chosen, the remainder of the solution can
be replaced by a regression. Formula (\ref{eq:final.main}), the first
order optimal prediction formula, is also a regression formula (an
average conditioned by partial information). To illustrate regression,
and more importantly, to remind the reader of formulas that will be
used in the sequel, we collect in the present section some regression
formulas for Gaussian measures, both for the discrete and the
continuous case. More details can be found in standard books
(e.g.\cite{papoulis} as well as in \cite{ckk3}).

We start by describing how to perform regressions on discrete sets of
Gaussian (normal) variables.  Let ${\mathbf u} = (u_1,\ldots,u_n)$ be
a real vector of jointly normal random variables; it has a probability
density $f({\mathbf u})$ of the form,
\begin{eqnarray}
\lefteqn{
P(s_1 < u_1 \le s_1 + ds_1, \dots, s_n < u_n \le s_n+ds_n) = 
f({\mathbf s}) \, ds_1\dots ds_n = } \nonumber \\
&& = Z^{-1} \exp \brk{- \half (\bfs,A\bfs) + {\bfb \cdot\bfs}} 
ds_1\dots ds_n,
\label{App_density}
\end{eqnarray}
where $Z$ is the appropriate normalization factor,
$\bfs=(s_1,\dots,s_n)$, and the $n\times n$ matrix $A$ with entries
$a_{ij}$ is symmetric, positive definite,
and has an inverse $A^{-1}$ . The matrix
$A^{-1}$ is the pairwise covariance matrix with elements
\[
a^{-1}_{ij} = \Cov{u_i}{u_j} \equiv \Average{u_i u_j} - \Average{u_i}
\Average{u_j},
\]
where the brackets, $\left< \cdot \right>$, denote averaging with
respect to the probability density; the vector ${\mathbf b}$ with
components $b_i$ is related to the expectation values of $\bfu$,
$\left<u\right>=(\left<u_1\right>,\dots,\left<u_N\right>)$, by
\[
A^{-1}\bfb = \Average{u}.
\]
The distribution is fully determined by the $n$ means and by the
$\half n(n+1)$ independent elements of the covariance matrix,
and therefore all the expectation value of any observable can 
be expressed in terms of these parameters.

Next, we assume that the random vector ${\mathbf u}$ satisfies a set
of conditions of the affine form,
\begin{equation}
g_{\alpha i} u_i = V_\alpha,   \qquad \alpha=1,\ldots,N<n,
\label{App_conditions}
\end{equation}
where the index $\alpha$ enumerates the conditions and summation over
repeated indices is assumed. Each equation in (\ref{App_conditions})
is a discrete analog of one of the equations in (\ref{eq:eq2}).  The
$N\times n$ matrix $G$, whose entries are $g_{\alpha i}$, determines
the full set of conditions.  To distinguish between the random
variables $(u_1,\ldots,u_n)$, and the collective variables
$(V_1,\dots,V_N)$, we enumerate the former by Roman indices and the
latter by Greek indices.

Our goal is to compute regressions, $E[\phi(u)|V]$, for various
functions $\phi$, i.e., conditional expectation values, or
equivalently, averages over the functions that satisfy the conditions.
We state three lemmas that will become handy below; for proofs, see
\cite{ckk3}.

\begin{lemma} \label{lemma1}
The conditional expectation of the variables $u_i$ is an affine
function of the conditioning data $V_{\alpha}$:
\begin{equation}
\Average{u_i}_{V} = q_{i \alpha} V_\alpha + c_i,
\label{constrainedMean}
\end{equation}
where the $n \times N$ matrix $Q$ whose entries are the $q_{i\alpha}$
and the $n$-vector ${\mathbf c}$ whose entries are the $c_i$ are given
by:
\begin{eqnarray}
Q &=& (A^{-1} G^{\dagger})(G A^{-1} G^{\dagger})^{-1}, \nonumber \\
{\mathbf c} &=& A^{-1}{\mathbf b} - 
(A^{-1} G^T)(G A^{-1} G^{\dagger})^{-1}(G A^{-1}{\mathbf b}), 
\end{eqnarray}
where the dagger denotes a transpose.
\end{lemma}

\begin{lemma}   
\label{lemma2}
The conditional covariance matrix has entries
\begin{eqnarray}
\Cov{u_i}{u_j}_{V} & = & 
\Average{u_i u_j}_{V} - \Average{u_i}_{V}
\Average{u_j}_{V} = \nonumber \\
& = & \Brk{A^{-1} - 
(A^{-1} G^{\dagger})(G A^{-1} G^{\dagger})^{-1}(G A^{-1})}_{ij},
\label{constrainedCovariance}
\end{eqnarray}
where the subscript $[]_{ij}$ denotes the $\{i j\}$ component of a
matrix.
\end{lemma}

\begin{lemma}   
\label{lemma3}
Wick's theorem holds for constrained expectations, namely,
\begin{eqnarray}
\lefteqn{
\Average{\prod_{p=1}^{P} 
\brk{u_{i_p} -\Average{u_{i_p}}_{V}}}_{V} = }\nonumber \\
&& = \left\{
\begin{array}{ll}
0    &   P \mbox{ odd} \\
\sum_{\mbox{perm}} \Cov{u_{i_1}}{u_{i_2}}_{V} \cdots
\Cov{u_{i_{P-1}}}{u_{i_P}}_{V} &   P \mbox{ even} \\
\end{array} \right.    \label{constrainedWick}
\end{eqnarray}
where the summation is over all possible pairings of the $P$
coordinates that are in the list.
\end{lemma}

Equation~(\ref{constrainedMean}) shows that conditioning data alter
expectation values linearly in the $V_\alpha$ and independently of
multiplicative factors in the covariances.
Equation~(\ref{constrainedCovariance}) shows that conditioned
covariances are determined by the matrix $G$ alone, without reference
to the $V_\alpha$.  Equation~(\ref{constrainedWick}) shows that the
conditioned Gaussian distribution, while not satisfying the
requirement that the covariance matrix be non-singular, retains a key
property of Gaussian densities.

In the applications below we shall use Gaussian variables
parameterized by a continuous variable $x$, i.e., Gaussian random
functions $u=u(x)$. Their means, $\langle u(x)\rangle$, and
covariances, $a^{-1}(x,y)=\langle u(x)u(y)\rangle-\langle
u(x)\rangle\langle u(y)\rangle$, will be defined for all $(x,y)$ in an
appropriate range rather than only for integer values of $(i,j)$.  The
matrix $A^{-1}$ becomes the integral operator whose kernel is a
function $a^{-1}$.  The kernel $a=a(x,y)$ of the operator $A$ inverse
to $A^{-1}$ is defined by
\[
\int a^{-1}(x,y) a(y,z) \, dy = \delta(x,z).
\]
The vectors with entries $g_{\alpha i}$, become functions
$g_{\alpha}(x)$, and the conditions (\ref{App_conditions})
become equations (\ref{eq:eq2}).  The regression formula,
(\ref{constrainedMean}), then changes into
\[
\Average{u(x)}_{V} = \Average{u(x)} +  c_{\beta}(x)
\Brk{V_\beta - \Average{\int g_\beta(y) u(y) \, dy}},
\]
where
\[
c_\beta(x) = \BRK{\int a^{-1}(x,y) g_\alpha(y) \, dy}
m^{-1}_{\alpha\beta},
\]
and the $m^{-1}_{\alpha\beta}$ are the entries of the matrix $M^{-1}$
whose inverse $M$ has entries
\[
m_{\alpha\beta} = \int \int g_\alpha(x) a^{-1}(x,y) g_\beta(y) \, dx \, dy.
\]
The formula for the constrained covariance can be obtained
from~(\ref{constrainedCovariance}) by replacing each $i$ by an $x$,
each $j$ by a $y$, and each summation over a Latin index by the
corresponding integration. Wick's theorem is still valid with the
appropriate changes in notation; note that the Greek indices, which
refer to the $N$ initial data, remain integers.


\section{An overview of invariant measures and Hamiltonian formalisms} 
\label{sec:measures}
\setcounter{equation}{0}
\setcounter{figure}{0}
\setcounter{table}{0}

Before proceeding with our numerical program, we summarize some
material on Hamiltonian systems, invariant measures in finite and
infinite dimensional systems, and the properties of certain
measures. This material can be found in books on quantum field theory
and related topics (see, e.g.,
\cite{goldstein,mandl,ramond}). More specific references
will be given below; we take a very elementary point of view.

A Hamiltonian system is described in terms of $n$ ``position''
variables, $q_i$, and their associated ``momenta'', $p_i$,
$i=1,\dots,n$; a Hamiltonian function ${\cal H}={\cal
H}(q_i,p_i)={\cal H}(q,p)$ is given, and the equations of motion are
\begin{equation}
\frac {dq_i}{dt}={\cal H}_{p_i},\quad \frac{dp_i}{dt}=-{\cal H}_{q_i},
\qquad i=1,\ldots,n.
\label{discH}
\end{equation}
If the initial values of the $2n$ variables $q,p$ are given, it is
assumed that the system (\ref{discH}) has a global solution in time.
Suppose the initial data are chosen at random in that $2n$ dimensional
space, with a probability density $f(q,p,0)$; it is easy to check that
the probability density of the $q$'s and $p$'s at time $t$ satisfies
the Liouville equation (see \cite{goldstein}):

\begin{equation}
f_t+\sum_{i}\left[\frac{dq_i}{dt}f_{q_i}+\frac{dp_i}{dt}f_{p_i}\right]=0.
\label{liouville}
\end{equation}

An invariant density is a probability density that does not depend on
time, i.e., one that satisfies equation (\ref{liouville}) with
$f_t=0$.  Recall that by the definition of a probability density, $f
\ge 0$, and $\int f \,dp \,dq=1$ (with obvious notations). One can
readily see that any function of ${\cal H}$ with these two properties
is an invariant density; the one that is natural for physical reasons
is $f=Z^{-1}\exp(-{\cal H}/T)$, where $Z$ is a normalization constant
and $T$ is the ``temperature''.  We set $T=1$ unless specified
otherwise.  If the initial data are sampled from this initial
distribution, and each of these samples is used as an initial datum
for the equations of motion, then the probability distribution of the
variables $q$ and $p$ at any later time $t$ will be the same as it was
initially. We now wish to generalize these notions of Hamiltonian
systems and invariant distributions to the infinite-dimensional case,
where the equations of motion will be partial differential equations
and the invariant distributions will be called ``invariant measures''.
We do so by way of an example that will be used in later sections.

Take the interval $[0,2\pi]$ and divide it into $n$ segments of length
$h$. At each mesh point $jh$, $j=1,\dots,n$ define variables $q_i,p_i$
; introduce the ``Hamiltonian''
\begin{equation}
{\cal H} = \sum_{i=1}^{n} h \Brk{
\frac{(p_{i+1}-p_i)^2}{2 h^2} +
\frac{(q_{i+1}-q_i)^2}{2 h^2} + F(p_i,q_i)},
\label{diffH}
\end{equation}
where values of $q,p$ outside the interval are provided by an
assumption of periodicity, and the term in brackets is a Hamiltonian
density. Consider the set of ordinary differential equations:

\begin{equation}
\frac{dp_i}{dt} = -\frac{1}{h}{\cal H}_{q_i} = \Delta_h q - F_{q_i},
\label{diff1}
\end{equation}
where $\Delta_h q = (q_{i+1}-2q_i+q_{i-1})/h^2$, and similarly,
\begin{equation}
\frac{dq_i}{dt} = \frac{1}{h}{\cal H}_{p_i} = -\Delta_h p + F_{p_i}.
\label{diff2}
\end{equation}
Note that the right hand side contains a factor $h^{-1}$ that has no analog
in the finite dimensional system above; its effect is to differentiate
the Hamiltonian density rather than the Hamiltonian itself.  This
modification is needed to get self-consistent limits as $h\to0$, a
limit operation we shall now undertake (see
\cite{goldstein,mandl}). As $h\to0$, these equations formally converge
to:
\begin{equation}
\left\{ 
\begin{array}{l}
\displaystyle{p_t =   q_{xx} - F_q,}  \\ \\
\displaystyle{q_t = - p_{xx} + F_p,}
\label{pde}
\end{array}
\right.
\end{equation}
or, writing $u=q+\imath p$ with imaginary $\imath$ and
$F'(u)=F_q+\imath F_p$, we find an equation of Schr\"odinger type:
\begin{equation}
\imath u_t = -u_{xx} + F'(u).
\label{pde2}
\end{equation}
The Hamiltonian, (\ref{diffH}), converges formally to 
\begin{equation}
\int_0^{2\pi}\Brk{\frac{1}{2} |u_x|^2 + F(u)} \,dx,
\label{formal}
\end{equation}
where the vertical lines denote a modulus.  One has to examine what in
these passages to a limit is justifiable, see e.g. \cite{ramond} for a
physics analysis, and \cite{mckean} for a mathematical analysis.

For every finite $h$, the (finite dimensional) measure
\begin{equation}
Z^{-1}\exp \Brk{-h \sum \brk{\frac{|u_{i+1}-u_i|^2}{2h^2} + F}}
\label{expdiff},
\end{equation}
is invariant for the system (\ref{diff1}--\ref{diff2}). Note that this
is true both with and without the extra factor $h^{-1}$ in equations
(\ref{diff1}--\ref{diff2}).  What is the limit of this measure as
$h\to0$?  Set for a moment $F=0$.  The exponential in (\ref{diffH})
factors into a product of terms each one of which is the exponential
of a single difference quotient, of the form
$\exp[-(p_{i+1}-p_i)^2/2h] $ or with $q$ replacing $p$.  Hence, the
$p$ variables are independent of the $q$ variables.  Furthermore, the
``increments'' $p_{i+1}-p_i$, $q_{i+1}-q_i$, are obviously Gaussianly
distributed, have a variance proportional to the distance $h$, and are
all independent of each other.  Thus in the limit, the functions
$q(x)$, (and similarly for $p(x)$), are made up of independent
Gaussian increments.  They differ from Brownian motion (see
\cite{freedman}) by being periodic rather than satisfying $q(0)=0$
(they are ``Brownian bridges'' -- which does not make a deep
difference).  Also, as long as $F=0$, the common value of $q(0)$ and
$q(2\pi)$ are undetermined because the exponent of the exponential is
unchanged when one adds a constant to the $q_i$; one can remove this
degeneracy by adding a term to the exponent that is sensitive to the
value of $q(0)$.  Thus, the limit of
\[
Z^{-1}\exp
\BRK{-h\sum\Brk{
\frac{(q_{i+1}-q_i)^2}{2h^2} +\frac{(p_{i+1}-p_i)^2}{2h^2} + F(q,p)}} \, 
dq_1 \cdots d q_n dp_1\cdots dp_n ,
\]
can be written as
\begin{equation}
dB_c\cdot \exp(-\int F dx),
\label{dBc}
\end{equation}
where $B_c$ is a suitably conditioned Brownian (Wiener) measure.
 
As is well-known (see, e.g., \cite{freedman}), a sample Brownian path
is, with probability one, nowhere differentiable. This fact is
not changed by the factor $\exp(-\int F dx)$ in (\ref{dBc}); thus if
we sample initial data from (\ref{dBc}) we obtain weak solutions of
the equation of motion (\ref{pde2}), as was indeed assumed in section
\ref{sec:weak}.
In addition, the integral $\int |u_x|^2dx$ diverges, so that the limit
in (\ref{diffH}) is purely formal; its meaning is given by equation
(\ref{dBc}). An important consequence of these facts is that the
exponential in (\ref{expdiff}) tends to zero as $h\to0$; this is
indeed necessary if we are to have a reasonable function-space
measure: As $h\to0$, we have a measure on a space of increasing
dimension, the density of functions that satisfy the set of
inequalities $s_i < q_i \le s_i + ds_i, i=1,\dots,n$, should decrease
as $n$ increases, thus $\exp(-{\cal H})$ should tend to zero and
${\cal H}$ should diverge.

Weak solutions that are spatially like Brownian motion are difficult
to resolve; difference quotients do no converge, and the Fourier
series expansions of these solutions converge slowly.  We are thus
consistent: our machinery applies in problems where it is indeed
needed.  Such problems are not exceptional; for example, in the
vanishing-viscosity limit, the solutions of the Euler equations have a
H\"older exponent of $1/3$, i.e., they are even less smooth than
Brownian motion (see \cite{onsager}).

Finally, a computational comment that will be useful below: Brownian
motions and Brownian bridges are easy to sample via interpolation
formulas related to the regression formulas of the previous Section
\ref{sec:regression} (see e.g. \cite{chorin},\cite{morokoff}); to modify
the measures so as to take into account the factor $\exp(-\int F dx)$
requires simply that the samples be rejected or accepted with some
form of a Metropolis algorithm.


\section{Conditional expectations with a non-Gaussian prior}    
\label{sec:non-gaussian}

\setcounter{equation}{0}
\setcounter{figure}{0}
\setcounter{table}{0}

Our goal in this section is to introduce a systematic approach for
solving the equations of optimal prediction (\ref{eq:final.main}) for
nonlinear equations of the form (\ref{eq:basic}).  We assume that
equation (\ref{eq:basic}) has Hamiltonian form, i.e., we assume that
there exists a Hamiltonian ${\cal H}$ such that (\ref{eq:basic}) is
the Hamilton equation of motion with this Hamiltonian.  We can
therefore assume the existence of a prior measure, $\mu(u)$ whose
density $f$ has the form
\begin{equation}
f(u) = Z^{-1} e^{-{\cal H}(u)},   
\label{eq:prior.measure}
\end{equation}
with $Z$ a normalization constant.  

In order to solve equation (\ref{eq:final.main}), we must compute the
conditional expectations on its RHS,
\begin{equation}
\brk{g,\Average{R(u)}_{V}},  
\label{eq:rhs}
\end{equation}
with $R$ being the RHS of equation (\ref{eq:basic}).  
This task is relatively simple when the prior measure $\mu$ is Gaussian; 
we address here the problem of what to do when it is not.  
The method we present is based on perturbation theory, 
see e.g. \cite{fetter-walecka:quantum,kleinert:gauge,ma}.  
The idea is to reduce the computation of (\ref{eq:rhs}) to the 
computation of regressions with respect to a Gaussian measure via 
a perturbation expansion.

We start by splitting the Hamiltonian into two parts,
\begin{equation}
{\cal H} = {\cal H}^{0} + {\cal H}^{1}.   
\label{eq:split.two}
\end{equation}
Here ${\cal H}^{0}$ is quadratic (producing a Gaussian measure
$\mu^0$) and ${\cal H}^{1}$ is a non-quadratic perturbation.  A
conditional expectation of a functional ${\cal F}(u)$ is defined as
\begin{equation}
\Average{{\cal F}}_{V} = Z^{-1} \int {\cal F}(u) d\mu_{V},
\label{eq:general.expect}
\end{equation}
with the normalization constant $Z = \int d \mu_{V}$.  Averages with
respect to Gaussian measures will be singled out by a superscript $0$:
\begin{equation}
\Average{{\cal F}}^0_V = Z_{0}^{-1} \int {\cal F}(u) d \mu_V^0,
\label{eq:gaussian.expect}
\end{equation}
with $Z_{0} = \int d \mu_V^0$.  Based on the Hamiltonian split,
(\ref{eq:split.two}), we can carry out the expansion
\begin{equation}
e^{-{\cal H}} = e^{-{\cal H}^{0}} \sum_{k=0}^{\infty}
\frac{(-1)^{k}}{k!}\brk{{\cal H}^{1}}^{k},  
\label{eq:hamiltonian.split}
\end{equation}
which can be then utilized to write
\begin{equation}
\int {\cal F}(u) d\mu_{V} = \sum_{k=0}^{\infty} \frac{(-1)^{k}}{k!} 
\int \Brk{{\cal F}(u)({\cal H}^{1})^k} d\mu_V^0.
\label{eq:last.int}
\end{equation}
Hence, combining (\ref{eq:general.expect}), (\ref{eq:gaussian.expect})
and (\ref{eq:last.int}) we find
\[
\Average{{\cal F}}_{V} = 
\sum_{k=0}^{\infty} \frac{(-1)^{k}}{k!} \frac{Z_{0}}{Z}
\Average{({\cal H}^{1})^{k} {\cal F}}_V^0.
\]
The ratio $Z_{0}/Z$ is:
\[
\brk{\frac{Z_{0}}{Z}}^{-1} =
\frac{\int d\mu_{V}}{\int d\mu_{V}^{0}} = 
\sum_{k=0}^{\infty} \frac{(-1)^{k}}{k!} 
\frac{\int({\cal H}^{1})^{k} d\mu_{V}^{0}}{\int d\mu_{V}^{0}}
= \sum_{k=0}^{\infty} \frac{(-1)^{k}}{k!} 
\Average{{\cal H}^{1})^{k}}_V^0,
\]
and therefore:
\[
\Average{{\cal F}}_{V} = 
\frac{\sum_{k=0}^{\infty} \frac{(-1)^{k}}{k!} 
\Average{({\cal H}^{1})^{k} {\cal F}}_V^0}
{\sum_{k=0}^{\infty} \frac{(-1)^{k}}{k!} 
\Average{({\cal H}^{1})^{k}}_V^0}.
\]
In particular, for $\alpha=1, \ldots,N$,
\begin{equation}
\brk{ g_{\alpha},\Average{R(u)}_{V}}  =  
\frac{\sum_{k=0}^{\infty} \frac{(-1)^{k}}{k!}
\brk{g_{\alpha}, \Average{R(u) \left({\cal H}^{1}\right)^{k}}_V^0}}  
{\sum_{k=0}^{\infty} \frac{(-1)^{k}}{k!}
\Average{({\cal H}^{1})^{k}}_V^0}.  
\label{eq:perturb}
\end{equation}
The RHS of (\ref{eq:perturb}) is now written in terms of expectation
values that we already know how to compute since they are averages
with respect to a Gaussian measure.  Note that the division by $Z/Z_0$
can be avoided by removing certain terms from the numerator of
(\ref{eq:perturb}); indeed, (\ref{eq:perturb}) is a conditioned
expansion in Feynman diagrams and it can be normalized by removing
unconnected graphs 
(see \cite{fetter-walecka:quantum,kleinert:gauge,ma}). 
We choose not to explain this fact here.  Note also that
the leading term in the expansion, (corresponding to $k=0$), where the
measure $\mu_V$ is simply replaced by $\mu_V^0$, already contains a
contribution of the nonlinear terms in the equation of motion.

Computing a finite number of terms in the series expansion
(\ref{eq:perturb}) can still be a relatively complicated task, and we
demonstrate in \S\ref{sec:example} a step-by-step solution of a model
problem.  The reader should not be unduly worried by the complexity of
some of the expressions, because : (i) To make the exposition as clear
as possible, no advantage is being taken here of various ways of
simplifying the expressions, such as, e.g., using orthogonal functions
for the $g_\alpha$; (ii) much of the algebra can be automated (see
below).

The partition of ${\cal H}$ into Gaussian and non-Gaussian parts is
not unique, and in section \ref{sec:example} and \ref{sec:dealiasing}
below we shall discuss
some ways to optimize it in order to gain accuracy.


\section{A perturbative treatment of a nonlinear Schr\"odinger
equation}
\label{sec:example}
\setcounter{equation}{0}
\setcounter{figure}{0}
\setcounter{table}{0}

We now utilize the perturbation method to approximating solutions of
the non-linear Schr\"odinger equation of the form (\ref{pde2}),
\begin{equation}
\imath u_t = - u_{xx} + \frac{1}{4}\Brk{3\Abs{u}^2 u + u^{*3}}, 
\quad  u=q+\imath p,
\label{eq:complex.schrodinger}
\end{equation}
in the interval $[0,2\pi]$, with periodic boundary conditions. 
Equation
(\ref{eq:complex.schrodinger}) can also be written as the pair of
equations,
\begin{equation}
\left\{
 \begin{array}{l}
\displaystyle{p_{t} =   q_{xx} - q^{3},}  \\ \\
\displaystyle{q_{t} = - p_{xx} + p^{3}.}   
 \end{array}
\right.
\label{eq:schrodinger}
\end{equation}
The corresponding Hamiltonian is
\begin{equation}
{\cal H}(p,q) = 
\frac{1}{2} \int \brk{p_{x}^{2}(x) + q_{x}^2(x) + 
\frac{1}{2}\Brk{p^4(x)+q^4(x)} }\, dx, 
\label{eq:sch.ham}
\end{equation}
and hence, equations (\ref{eq:schrodinger}) preserve the canonical
density
\[
f_{0}(p,q) = Z^{-1} e^{-{\cal H}(p,q)}.
\]

(Note that the temperature $T$ has been set equal to 1).  To simplify
the example, we follow \cite{ckk3} and use the same kernels in the
definition of the collective variables for $p$ and $q$, i.e., the
collective variables are
\begin{equation}
\{U_{\alpha}^{p},U_{\alpha}^{q}\} = 
\BRK{(g_{\alpha},p),(g_{\alpha},q)}, \qquad \alpha =1,\ldots,N,   
\label{eq:sch.collective}
\end{equation}
and their initial values, $V_{\alpha}^p,V_{\alpha}^q$, are given.  The
system of ordinary differential equations arising out of equation
(\ref{eq:final.main}) for the $V_{\alpha}^p,V_{\alpha}^q$ is:
\begin{equation}
\left\{ 
\begin{array}{l}
\displaystyle{\frac{d V_{\alpha}^{p}}{dt} =
+\brk{g_{\alpha}, \frac{\partial^{2}}{\partial x^{2}} \Average{q}_{V}} - 
\brk{g_{\alpha}, \Average{q^{3}}_{V}},} \\ \\
\displaystyle{\frac{d V_{\alpha}^{q}}{dt} = 
-\brk{g_{\alpha}, \frac{\partial^{2}}{\partial x^{2}} \Average{p}_{V}}
+ \brk{g_{\alpha}, \Average{p^{3}}_{V}}.}
\end{array} \right. 
\qquad \alpha = 1,\ldots,N.  
\label{eq:sch.ode}
\end{equation}
We therefore have to compute the four terms 
$\brk{g_{\alpha}, \partial_{xx} \Average{q}_{V}}$,
$\brk{g_{\alpha}, \partial_{xx} \Average{p}_{V}}$,
$\brk{g_{\alpha}, \Average{q^{3}}_{V}}$ and
$\brk{g_{\alpha}, \Average{p^{3}}_{V}}$.

The first step is to split the Hamiltonian into two parts - a
quadratic and a non-quadratic part, ${\cal H} = {\cal H}^{0} + {\cal
H}^{1}$; This will be done in the next section.  We will just note at
that point, that ${\cal H}^{0}$ will be of the form ${\cal H}^{0} =
\frac{1}{2} \int \left[ p_{x}^{2}(x) + q_{x}^2(x) + m_{0}^2
\left(p^2(x)+q^2(x)\right) \right] dx$.  Once this is done, The RHS of
(\ref{eq:sch.ode}) can be computed following (\ref{eq:perturb}).  For
example, the term $(g_{\alpha}, \Average{p^{3}}_{V})$
can be obtained by substituting $p^3$ for $R(u)$ in that equation. 
In particular, 
the zeroth, leading term is:
\begin{equation}
\brk{g_{\alpha}, \Average{p^{3}}_V^0},  
\label{eq:sch.zero}
\end{equation}
(note that this term already includes a nonlinear effect); a
first-order (in the perturbation series) approximation will add the
term
\begin{equation}
\frac{
\brk{g_{\alpha}, \Average{p^{3} 
{\cal H}^{1}}_V^0}}{\brk{1-\Average{{\cal H}^{1}}_V^0}},
\label{eq:sch.one}
\end{equation}
and so on, with similar expressions for the rest of the terms on the
RHS of (\ref{eq:sch.ode}).  Note that all these expansions are used to
solve the equations of first-order prediction; in principle, we can
perform higher-order predictions by using higher moments in setting up
the matching of conditions before equation (\ref{eq:final.main}), and
then improve the evaluation of the right-hand side of equation
(\ref{eq:final.main}) be using more terms in the perturbation
expansion.  We consider here only the latter possibility.

The leading-order terms (\ref{eq:sch.zero}) can be computed
using the results of \S\ref{sec:regression}.  For the first-order
term, (\ref{eq:sch.one}), we have
\[
 \left<p^{3}{\cal H}^{1} \right>_V^0 = 
 \frac{1}{4} \left< p^{3}\int 
\left[p^{4}+q^{4}-2 m_{0}^{2}(p^{2}+q^{2}) \right]dz\right>_V^0;
\]
and we therefore have to compute
\[
 \left<p^{3}\int p^{s}dz \right>_V^0 = 
\int_{z} \left< p^{3}(x) p^{s}(z) \right>_V^0 dz, \qquad
 \left<p^{3}\int q^{s}dz \right>_V^0 = 
\int_{z} \left< p^{3}(x) q^{s}(z) \right>_V^0 dz
\]
for $s=2,4$.  To summarize, the RHS of (\ref{eq:sch.ode}), involves
integration (with respect to $z$) and differentiation (with respect to
$x$) of the following (with $j=1,3$ and $s=2,4$)
\begin{equation}
  \Average{p^{j}(x)p^{s}(z)}_V^0, \quad  
  \Average{p^{j}(x)q^{s}(z)}_V^0, \quad
  \Average{q^{j}(x)p^{s}(z)}_V^0, \quad  
  \Average{q^{j}(x)q^{s}(z)}_V^0.  
  \label{eq:sch.moments}
\end{equation}
The terms in (\ref{eq:sch.moments}) can be computed by application of
Wick's theorem (Lemma \ref{lemma3}); the algebra can be performed with
aid of a symbolic computer program (such as Mathematica). In
particular, one obtains the identities:
\bea
\lefteqn{\left< p(x)\,p^{2}(z) \right> = 
 -2\left<p(x)\right>\left<p(z)\right>^2 
 + 2\left<p(z)\right>\left<p(x)p(z)\right> 
 + \left<p(x)\right>\left<p^{2}(z)\right>.} \nonumber \\
\;
\lefteqn{\left< p(x)\,p^{4}(z) \right> = 
6\left<p(x)\right>\left<p(z)\right>^4 - 
8\left<p(z)\right>^3\left<p(x)p(z)\right> -} \nonumber \\
&& - 12\left<p(x)\right>\left<p(z)\right>^2\left<p^{2}(z)\right> 
   + 12\left<p(z)\right>\left<p(x)p(z)\right>\left<p^{2}(z)\right> 
   + 3\left<p(x)\right>\left<p^{2}(z)\right>^2. \nonumber \\
\;
\lefteqn{\left<p^{3}(x)\,p^{2}(z) \right> = 
 -12\left<p(x)\right>^2\left<p(z)\right>\left<p(x)p(z)\right> + 
  6\left<p^{2}(x)\right>\left<p(z)\right>\left<p(x)p(z)\right> + } \nonumber \\
&& + \left<p(x)\right>^3\bigg(6\left<p(z)\right>^2 - 2\left<p^{2}(z)\right>\bigg) + \nonumber \\
&& + \left<p(x)\right>\left[6\left<p(x)p(z)\right>^2 + 
  \left<p^{2}(x)\right>\bigg(-6\left<p(z)\right>^2 + 
  3\left<p^{2}(z)\right>\bigg) \right]. \nonumber \\
\;
\lefteqn{\left< p^{3}(x)\,p^{4}(z) \right> =
   72\left<p(x)\right>^2\left<p(z)\right>\left<p(x)p(z)\right>
    \bigg(\left<p(z)\right>^2 - \left<p^{2}(z)\right>\bigg) + } \nonumber \\
&& + \left<p(x)\right>^3\bigg(-20\left<p(z)\right>^4 + 36\left<p(z)\right>^2
    \left<p^{2}(z)\right> - 6\left<p^{2}(z)\right>^2\bigg) + \nonumber \\
&& + 12\left<p(z)\right>\left<p(x)p(z)\right>\left[2\left<p(x)p(z)\right>^2 + 
    \left<p^{2}(x)\right>\bigg(-2\left<p(z)\right>^2 + 
    3\left<p^{2}(z)\right>\bigg)\right]. 
    + \nonumber \\
\;
&& + 9\left<p(x)\right>\left[4\left<p(x)p(z)\right>^2
    \bigg(-2\left<p(z)\right>^2 + \left<p^{2}(z)\right>\bigg) + \right. \nonumber \\
&& \left. + \left<p^{2}(x)\right>\bigg(2\left<p(z)\right>^4 
  - 4\left<p(z)\right>^2\left<p^{2}(z)\right> + \left<p^{2}(z)\right>^2\bigg)\right].
 \label{wick:expressions}
\eea
All the averages in equation (\ref{wick:expressions}) are of course
conditioned by $V$; the constant repetition of the subscript $V$ and
the superscript $0$ has been avoided for esthetic reasons.

We now pick the kernels $g_{\alpha}(x)$ to be translates of a fixed
function $g(x)$, i.e., $g_{\alpha}(x) = g(x-x_{\alpha})$, which is a
normalized (not random!) Gaussian with periodic images and width
$\sigma$:
\begin{equation}
 g(x) = \frac{1}{\sqrt{\pi}\sigma} \sum_{\tau = -\infty}^{\infty} 
 \exp \left[ - \frac{(x-2\pi\tau)^2}{\sigma^2} \right].
 \label{eq:cont.kernels}
\end{equation}
We note that the Fourier representation of $g(x)$ is
 $g(x)=\frac{1}{2\pi}\sum_{k=-\infty}^{\infty} e^{\imath k
 x-\frac{1}{4}k^2 \sigma^2}$.  Given this choice of kernel functions
 we can write
\begin{equation}
\left< p(x) \right>_V^0 = \sum_{\alpha=1}^{N} 
c_{\alpha}^{pp}(x) V_{\alpha}^{p},
\label{eq:p}
\end{equation}
where
\[
 c_{\alpha}^{pp}(x) = c_{\alpha}^{qq}(x) = 
  \frac{1}{2 \pi} \sum_{\beta=1}^{N} \sum_{k = -\infty}^{\infty}
  \frac{e^{-\frac{1}{4}k^2\sigma^2}}{k^2+m_{0}^{2}} 
  \exp[\imath k(x-x_{\beta})][m^{-1}]^{pp}_{\beta \alpha},
\]
and
\[
 m_{\beta \alpha}^{pp} = \frac{1}{2\pi} \sum_{k=-\infty}^{\infty}
 \frac{e^{-\frac{1}{2}k^2\sigma^2}}{k^2+m_{0}^{2}} \exp[\imath
 k(x_{\alpha}-x_{\beta})].
\]
This choice of kernels is the same as in previous work
\cite{ckk2,ckk3}. It is far from optimal in the context of perturbation theory;
in particular, orthogonal kernels such as the Fourier kernels used in
section 7 below reduce the number of non-zero terms in the
expansion. We thought that we should present at least once a
perturbative calculation with a general kernel.

\subsection{The partition of the Hamiltonian}   \label{subsection:partition}

We now turn to the question of how exactly the Hamiltonian should be
divided into a sum of a quadratic part and a perturbation, ${\cal H} =
{\cal H}^{0} + {\cal H}^{1}$. Of course we wish ${\cal H}^1$ to be as
small as possible, so as to have a perturbation series that behaves as
well as possible; we therefore perform a partition with a few free
parameters over which we shall minimize ${\cal H}^1$; we choose to
write:
\bea {\cal H}^0 &=&\frac{1}{2} \int_{0}^{2\pi} \left( |u_x|^2 +
m_0^2 \Abs{u}^2 + b\right) dx, \label{partition} \\ {\cal H}^1
&=&\frac{1}{2}
\int_{0}^{2\pi} \left( {F- m_0^2 \Abs{u}^2-b} \right) dx, \nonumber
\eea 
where, as before, $u=q+\imath p, F=\frac{1}{2}(p^4+q^4)$.  There
are no odd powers in the partition because the measure is invariant
under the reflection $u \leftrightarrow (-u)$; note that the term
$m_0^2 \Abs{u}^2$ removes the indeterminacy in the
Gaussian measure defined by ${\cal H}^{0}$. This is not the only
partition that can be considered; one could for example add and
subtract squares of fractional derivatives of $u$. The task at hand is
to choose good values for the parameters $m_0$ and $b$.  They cannot
in general be chosen so as to make the perturbation series convergent;
what one would really want is to make the measure defined by ${\cal
H}$ and conditioned by the $V_\alpha$ be a small perturbation of the
conditional measure defined by ${\cal H}^{0}$, and while this can in
principle be done, and would presumably lead to time-dependent
equations for $m_0$ and $b$, it is reasonable, as a first try, to
choose $m_0$ and $b$ so as to minimize $\left<{\cal H}^{1}\right>,
\left<({\cal H}^{1})^{2}\right>$, where the averages are
unconditional.  Note also that the presence of the term $b$ leaves all
averages with respect to the measure defined by ${\cal H}^{0}$
unchanged---it gets absorbed into the normalization constant $Z$---but
it does affect the expansion of $\exp(-{\cal H}^{1})$.

A straightforward algebraic manipulation, simplified by the fact that
$p$ and $q$ are uncorrelated and that $\left< p^{s} \right> = \left<
q^{s}
\right>$ for all $s$, leads to
\bea
\lefteqn{ \left< ({\cal H}^{1})^{2} \right>^{0} =
 \frac{1}{16} \int_{z=0}^{2\pi} \int_{x=0}^{2\pi} \,dx \,dz
\left[
2 \Average{p^{4}(x) p^{4}(z)}^0 + 
2 \brk{\Average{p^{4}(x)}^0}^2 - \right. } \nonumber \\
&& - 8 m_{0}^{2} \Average{p^{4}(x) p^{2}(z)}^0 -
8 m_{0}^{2} \Average{p^{4}(x)}^0 \Average{p^{2}(z)}^0 + \nonumber \\
&& \left. + 8 m_{0}^{4} \brk{\Average{p^{2}(x)}^0}^2 + 
8 m_{0}^{4} \Average{p^{2}(x)p^{2}(z)}^0 \right],  
\label{eq:sch.h1_2}
\eea
where the $\Average{\cdot}^0$, denotes average with respect to
the unconditional Gaussian measure.  Using the expressions
(\ref{wick:expressions}), one can numerically minimize
(\ref{eq:sch.h1_2}) as a function of $m_{0}$.  This minimum is
obtained at $m_{0} \approx 1.055.$

After setting $m_0$, we are free to pick the second constant, $b$, in
(\ref{partition}).  Changing $b$ will not affect the variance which we
just computed in (\ref{eq:sch.h1_2}).  We choose to set the first
moment of the perturbation to zero, i.e.,
\begin{equation}
\Average{{\cal H}^{1}}^0 = 0.  
\label{first.moment}
\end{equation}
Given the partition, (\ref{partition}), it is clear that
(\ref{first.moment}) holds if
\begin{equation}
b \pi = \frac{1}{4} \Average{\int \Brk{ p^{4}+q^{4} - 
2 m_{0}^{2}(p^{2}+q^{2}) } \,dx }^0.  
\label{b.pi}
\end{equation}
With the aid of Wick's theorem and equation (\ref{eq:p}), equation
(\ref{b.pi}) can be rewritten as
\begin{equation}
b \pi = \Brk{ \frac{3}{4\pi}\sum_{k=-\infty}^{\infty}
\frac{1}{k^{2}+m_{0}^2} - m_{0}^{2} } \cdot
\sum_{k=-\infty}^{\infty} \frac{1}{k^{2}+m_{0}^2}.
\label{eq:h1.tilde}
\end{equation}
If one chooses $m_{0} = 1.055$ so as to minimize the second moment of
$\left< {\cal H}^{1} \right>$, the RHS of (\ref{eq:h1.tilde}) equals
$-1.197$, which in turn determines $b=-0.381$ so that the first moment
of the perturbation will vanish.  Note that even the zero-th order
expansion uses information about the higher-order terms, since the
parameters that determine the partition and thus the zero-th order
term depend on an analysis of the later terms.

One should also note that once the future conditions have been determined,
one has to perform further regressions to obtain the mean solutions at
various points in space; the machinery there is exactly analogous to
what has just been done, and will not be spelled out here.

\subsection{Numerical checks}   \label{subsection:numerical.results}

We now check the perturbation series as well as the optimal prediction
scheme by comparing the results they give with numerical results
obtained at substantial expense by sampling the initial conditions,
solving the differential equations over and over, and averaging. We
concentrate in the present section on the case $N=2$, i.e., a case
where we have initially as data only four values of collective
variables, and are trying to find the mean future conditioned by these
four values.  We display only the variation in time of these
collective variables, $U_\alpha^{p}, U_\alpha^{q}, \alpha=1,2$, which
were defined in equation (\ref{eq:sch.collective}); The kernels are
taken as $g_{\alpha}(x) = g(x-x_{\alpha})$ with $x_{1} = \pi/8$ and
$x_{2} = 9\pi/8$.  The parameter $\sigma$ is set as $\pi$.  We first
use the perturbation series with the optimal values $m_0 = 1.055$ and
$b=-0.381$ computed in the previous section
\ref{subsection:partition}.  
We then display results obtained with different values of $m_0$ and
$b$.  Specifically, we choose $m_{0}=0.9$ and $b=0$, obtained by
splitting the Hamiltonian into
\bean {\cal H}^0 &=&\frac{1}{2} \int_{0}^{2\pi} \left( |u_x|^2 +
m_0^2 \Abs{u}^2 \right) dx, \\ 
{\cal H}^1 &=&\frac{1}{2}
\int_{0}^{2\pi} \left( {F- m_0^2 \Abs{u}^2} \right) dx, 
\eean
and requiring only that the first moment of the perturbation 
vanish, i.e., $\Average{{\cal H}^{1}}^0 = 0$ (compare 
with (\ref{partition}) and (\ref{first.moment})).

We check these results by replacing the continuum equations
(\ref{eq:schrodinger}), (\ref{eq:sch.ode}), by a formal finite
difference approximation conditioned by discrete forms of these
equations, and then display the convergence of the conditioned mean of
many solutions of the difference equations to the optimal prediction
obtained by the perturbation analysis.  Specifically, we replace
equation (\ref{eq:schrodinger}) by the following difference equations,
\begin{equation}
\left\{
\begin{array}{l}
\displaystyle{\frac{dp(j)}{dt} =  \frac{q(j-1)-2q(j)+q(j+1)}{h^{2}} 
  - q^{3}(j),} \\ \\
\displaystyle{\frac{dq(j)}{dt} =  - \frac{p(j-1)-2p(j)+p(j+1)}{h^{2}} 
  + p^{3}(j),} 
\end{array} \right. 
\qquad j=1,\ldots,n,  \label{eq:discrete.equations} 
\end{equation}
where $h=2 \pi / n$ is the mesh size.  
The conditions (\ref{eq:sch.collective}) are replaced
by the following discrete approximations: 
\begin{equation}
 U_{\alpha}^{p} = \sum_{j=1}^{n} h g_{\alpha}(j)p(j), \quad
 U_{\alpha}^{q} = \sum_{j=1}^{n} h g_{\alpha}(j)q(j), \qquad \alpha = 1,2,
\end{equation}

(A factor $h$ has been introduced in the definition of the collective
variables to allow them to converge to the continuum collective
variables $(g_{\alpha},u)=\int g_{\alpha}(x) u(x)dx$.)

The Hamiltonian (\ref{eq:sch.ham}) is replaced by the discrete
Hamiltonian:
\bea
 \lefteqn{{\cal H}(p,q) = } \label{eq:discrete.hamiltonian}  \\  
 && = \frac{h}{2} \sum_{j=1}^{n} 
  \left\{ \left[ \frac{p(j+1)-p(j)}{h} \right]^{2} + 
   \left[ \frac{q(j+1)-q(j)}{h} \right]^{2} +
	\frac{1}{2} [ p^{4}(j) + q^{4}(j)] \right\}  \nonumber
\eea

We present results obtained with two mesh sizes, corresponding to
$n=8$ and $n=16$.  We also checked that the results with $n=32$ are
very close to the results with $n=16$.  For each mesh size $h$ we use
a Metropolis Monte-Carlo algorithm to find 5000 initial data drawn
from the distribution defined by (\ref{eq:discrete.hamiltonian})
conditioned by the values of the collective variables, integrate the
equations in time up to $t=1$, and average the results at various
fixed time intervals.  This numerical calculation is very costly, even
for moderate values of $n$, but it is independent both of the
perturbative analysis and of the machinery of optimal prediction.  It
is important to note that as the mesh size $h$ tends to zero, the
results of a standard (i.e., non-averaged) finite-difference solution
of equation (\ref{eq:schrodinger}) with data drawn from the
distribution (\ref{eq:discrete.hamiltonian}) diverge pointwise as $h
\rightarrow 0$.  Conditional averaging provides the only meaningful
numerical solution of equations (\ref{eq:discrete.equations}) for such
initial data.

\begin{figure}
 \centerline{\includegraphics{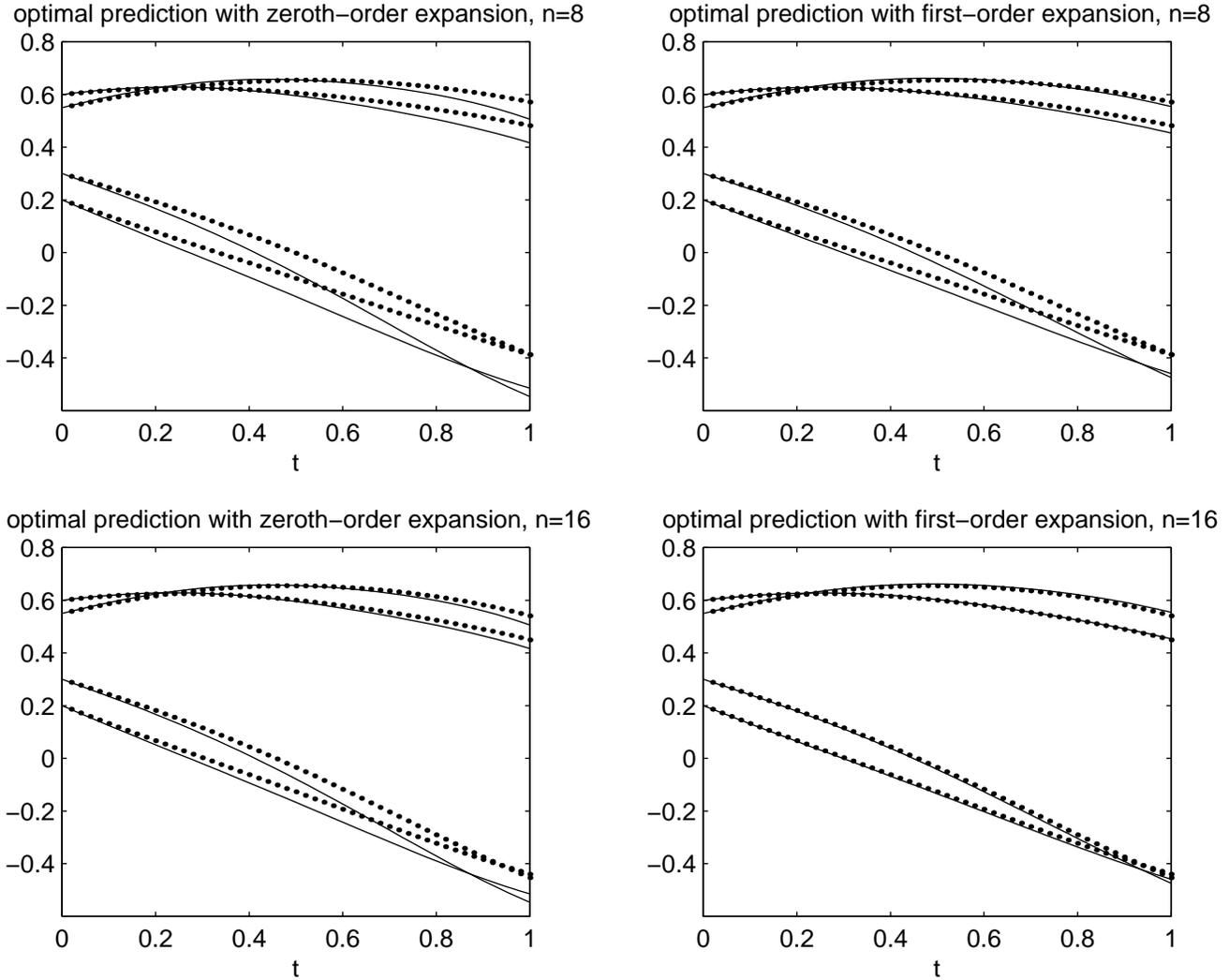}}
\caption{{\sf Time evolution of four collective variables 
 $U_{1}^{p}$, $U_{2}^{p}$, $U_{1}^{q}$ and $U_{2}^{q}$ for the
 nonlinear Schr\"odinger equation
 (\ref{eq:schrodinger}),(\ref{eq:discrete.equations}), with the
 optimal $m_0 = 1.055$, $b=-0.38$; Solid lines -- optimal prediction
 equations.  Dotted lines -- average over $5000$ solutions obtained
 from initial data sampled from the discrete Hamiltonian
 (\ref{eq:discrete.hamiltonian}) with $n=8$ and $n=16$
 points}. \label{figure:optimal}}
\end{figure}

\begin{figure}
 \centerline{\includegraphics{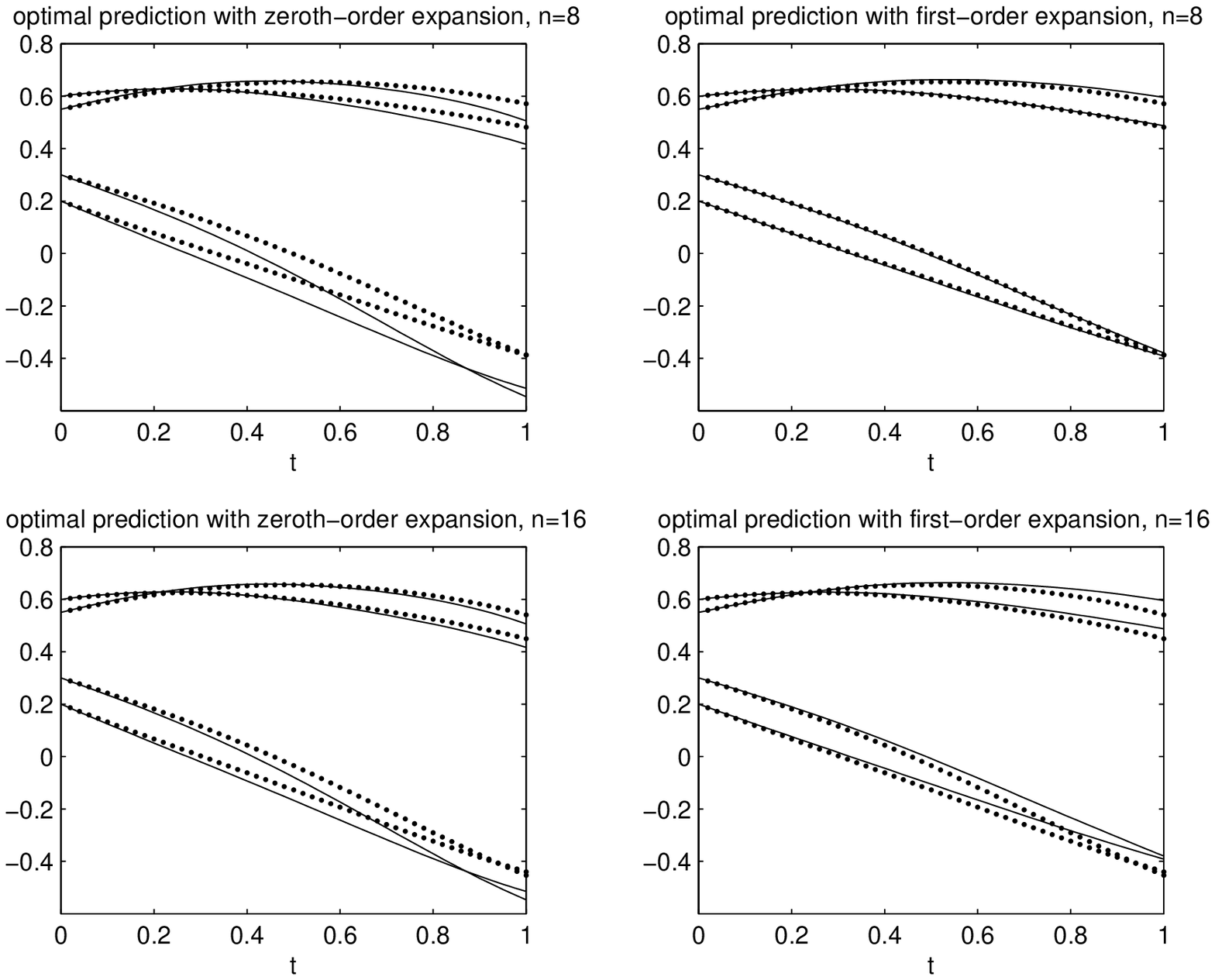}}
\caption{{\sf Time evolution of four collective variables 
 $U_{1}^{p}$, $U_{2}^{p}$, $U_{1}^{q}$ and $U_{2}^{q}$ for the
 nonlinear Schr\"odinger equation
 (\ref{eq:schrodinger}),(\ref{eq:discrete.equations}), with $m_0 =
 0.9$, $b=0$; Solid lines -- optimal prediction equations.  Dotted
 lines -- average over $5000$ solutions obtained from initial data
 sampled from the discrete Hamiltonian (\ref{eq:discrete.hamiltonian})
 with $n=8$ and $n=16$ points}. \label{figure:nonoptimal}}
\end{figure}

In Figure \ref{figure:optimal} we present the evolution in time of the
four collective variables $U_{1}^{p}$, $U_{2}^{p}$, $U_{1}^{q}$ and 
 $U_{2}^{q}$ with the optimal $m_0 = 1.055$, $b=-0.38$.  Figure
\ref{figure:nonoptimal} presents the plots corresponding to the
choice $m_0=0.9$ and $b=0$.  The zeroth-order solution is the optimal
prediction solution obtained with only the zeroth, leading, term in
the perturbation expansion (see, e.g., equation (\ref{eq:sch.zero})).
The first-order solution is the optimal prediction solution obtained
after adding a first-order correction to the perturbation series.
Note however that the optimal choice of parameters uses information
about terms of order one.

In both figures, the solid lines are
the solutions obtained with the optimal prediction equations, while the
dotted lines represent an average over
over $5000$ solutions that evolve from initial data sampled from 
the conditioned discrete Hamiltonian (\ref{eq:discrete.hamiltonian}).

Clearly, there is an improvement when one uses additional terms in the
perturbation series.  Also, even though the individual numerical
solutions converge only weakly to a continuum limit, the average over
numerical solutions with data sampled from the discrete Hamiltonian on
one hand and the solution of the optimal prediction equations on the
other hand get close as $n$ increases.  The perturbation expansion
converges rapidly; the key graph is the one on the bottom right of
Figure \ref{figure:optimal}: the comparison between the expansion up
to first-order with the average numerical solution with $n=16$.  One
conclusion we draw from these graphs and use in section 7 below is
that with the optimal partition of the Hamiltonian one can obtain an
accurate solution with only the zero-th term in the expansion; the
computation of the optimal parameters in the expansion uses the
first-order term.

\begin{figure}
\centerline{\includegraphics[scale=0.9]{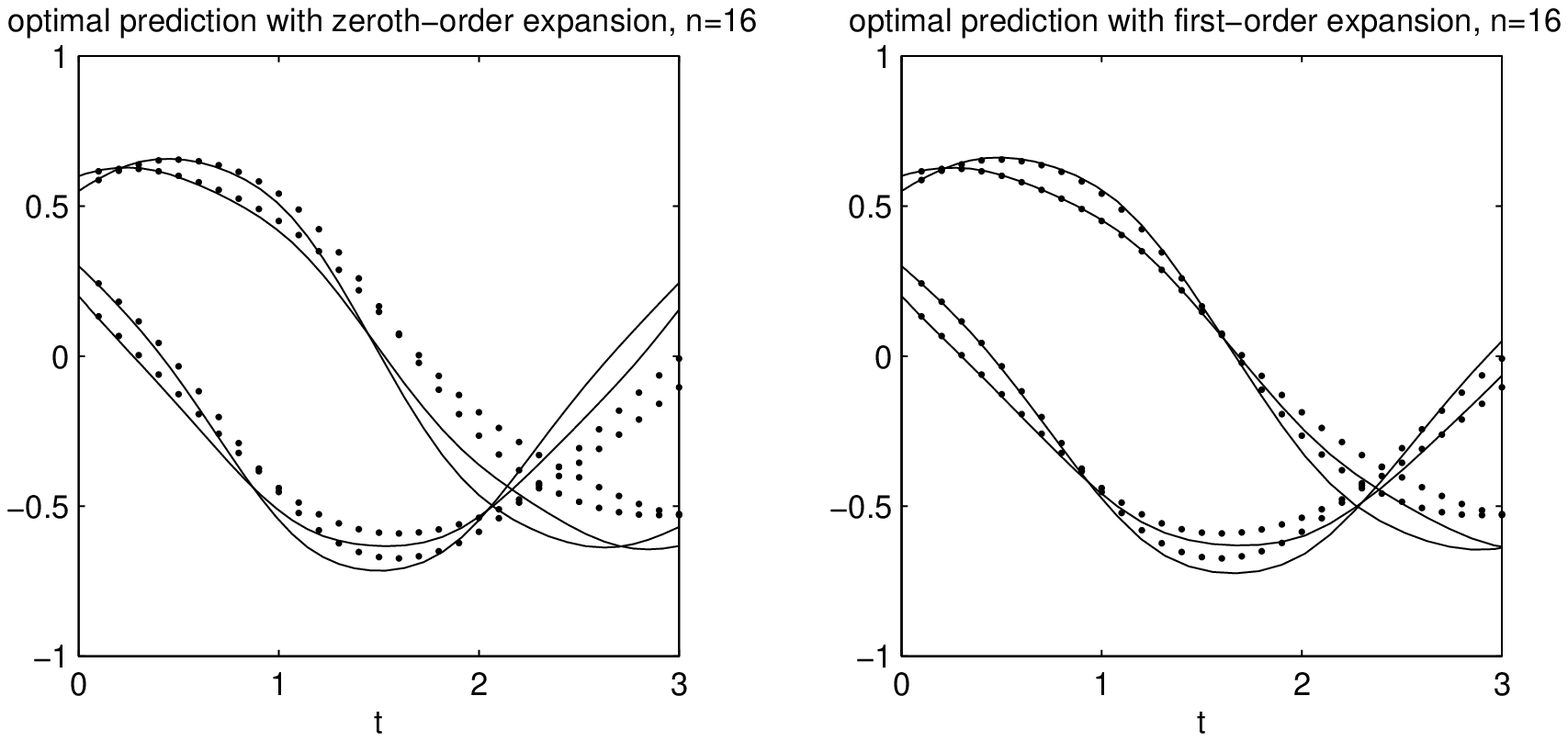}}
\caption{{\sf Longer time evolution of four collective variables 
$U_{1}^{p}$, $U_{2}^{p}$, $U_{1}^{q}$ and $U_{2}^{q}$ for the
nonlinear Schr\"odinger equation
(\ref{eq:schrodinger}),(\ref{eq:discrete.equations}), with the optimal
$m_0 = 1.055$, $b=-0.38$; Solid lines -- optimal prediction equations.
Dotted lines -- average over $5000$ solutions obtained from initial
data sampled from the discrete Hamiltonian
(\ref{eq:discrete.hamiltonian}) with $n=16$ points.}}
\label{figure:long}
\end{figure}

The limitations of first-order optimal prediction are displayed in
Figure \ref{figure:long}, where the integrations are carried out to
longer times.  We are working with a temperature $T=1$, i.e., the
fluctuation in $u$ are of order 1; by contrast, in the longer runs of
\cite{ckk2},\cite{ckk3} we used a smaller temperature $T=\pi/15$; we
also have only 4 collective variables, not enough at this temperature
to keep the variances of the collective variables small.  In Figure
\ref{figure:long} the optimal prediction solution is
based on the optimal choice of $m_0 = 1.055$, $b=-0.38$, and the
discrete solution is presented for $n=16$.  Once again, the numerical
solution is an average over 5000 individual solutions.  

\begin{figure}
\centerline{\includegraphics[scale=0.4]{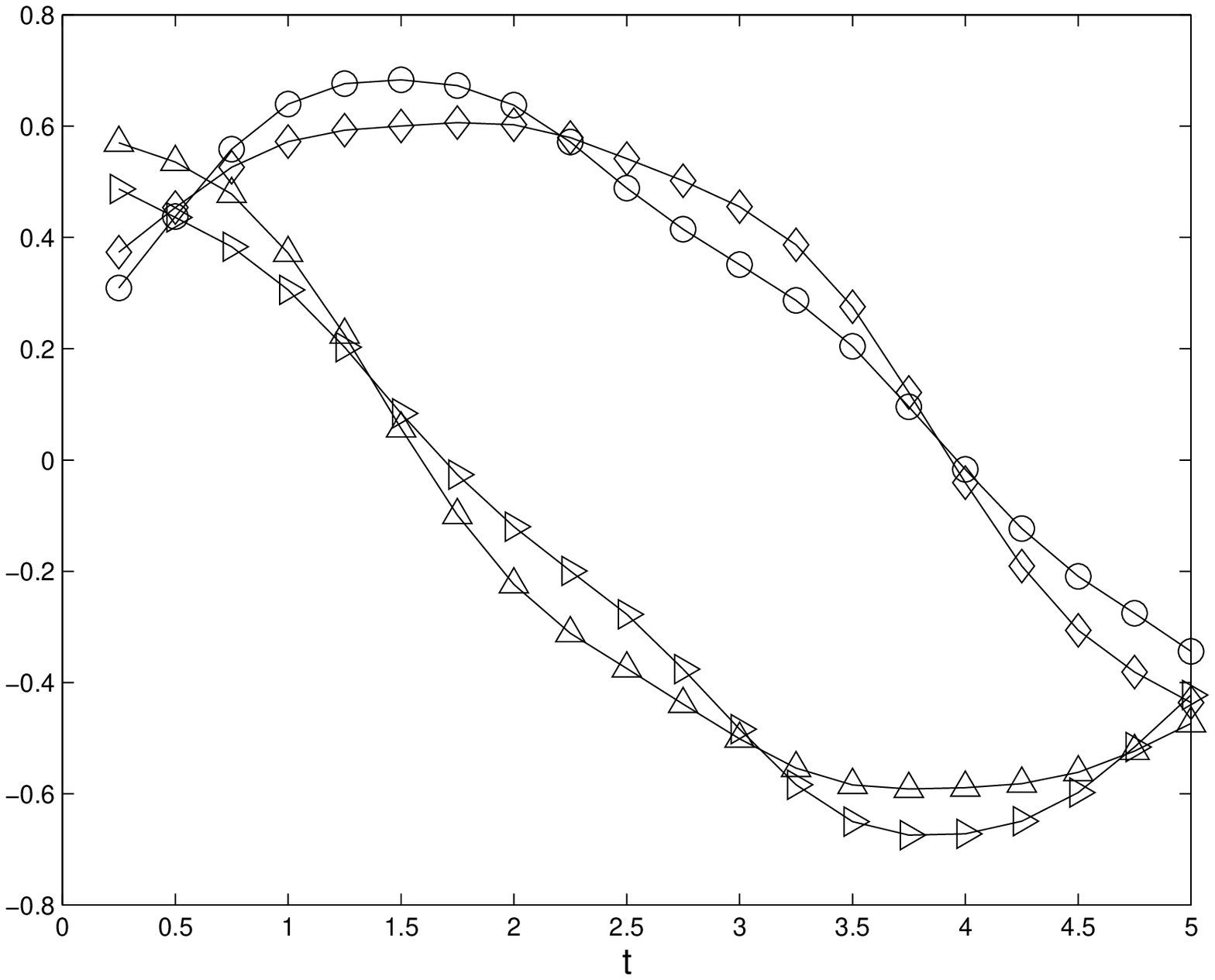}
            \includegraphics[scale=0.4]{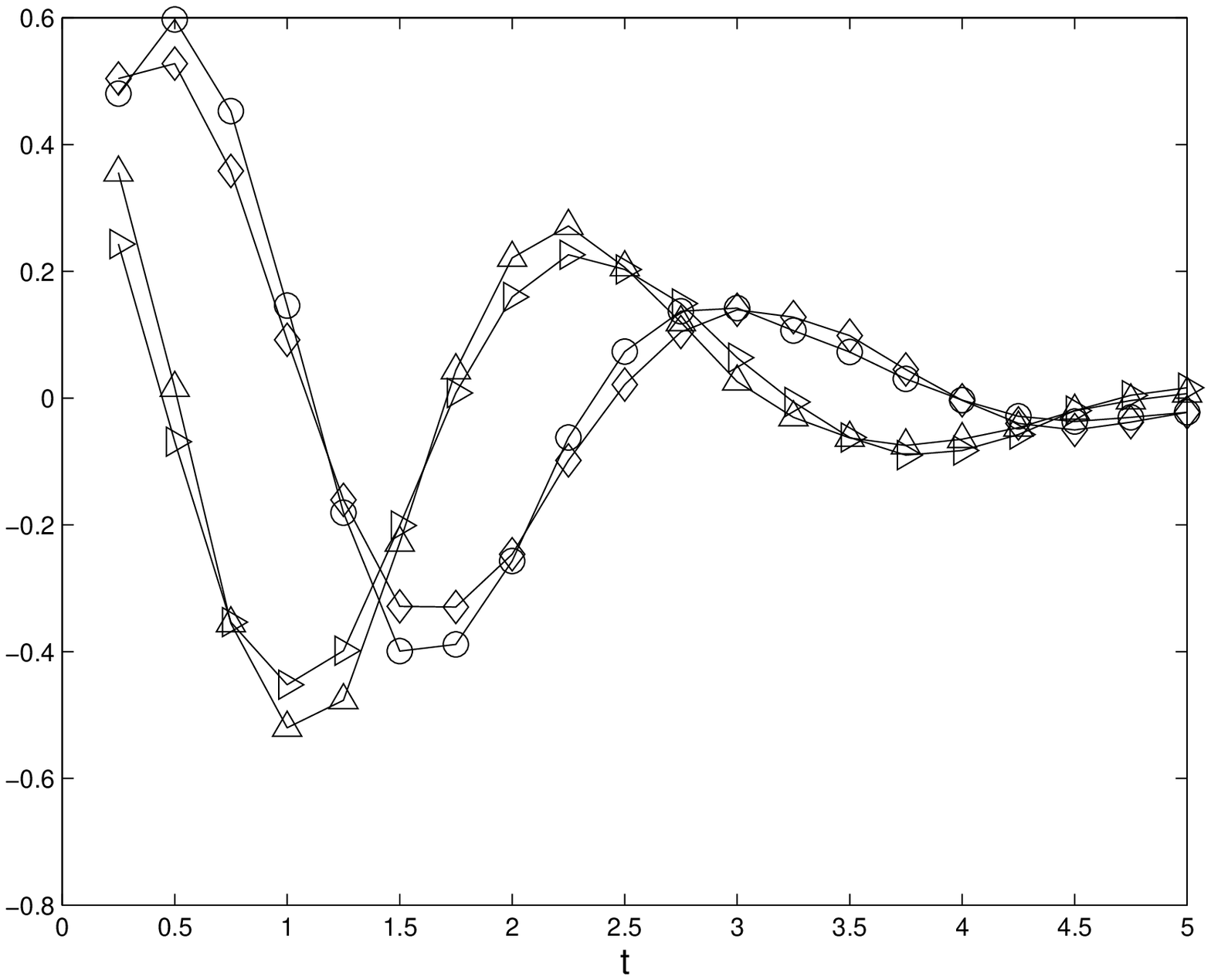}}
\caption{{\sf 
Long time evolution of the collective variables for initial
distributions of the form $e^{-{\cal H}/T}$ with (a) $T=0.2$ and (b)
$T=4$.  }}
\label{figure:decay}
\end{figure}

The effect of temperature is displayed in Figure \ref{figure:decay},
where we present the time evolution of the four collective variables,
$U^p_\alpha$, $U^q_\alpha$ up to time $t=5$. Both graphs are for the
same initial values of $U^p_\alpha$ and $U^q_\alpha$ but differ in the
temperature which determines the distribution of initial data. The
graph on the left corresponds to low temperature, $T=0.2$, whereas the
graph on the right corresponds to high temperature, $T=4$. As
expected, the higher the temperature, the faster is the decay towards
equilibrium; the collective variables tend faster towards their
equilibrium value of zero, and the first-order prediction scheme with a
small, fixed number of conditions loses accuracy faster. 
There are two ways to improve the prediction: Go to more
sophisticated prediction theory, as outlined in section 2, or increase
the number of collective variables. The first alternative will be
explored in later publications; the value of the second approach will
be shown in the next section, with a choice of kernels that reduces
the amount of labor and also makes possible an analytical estimate of
the difference between optimal prediction and a simple scheme that
makes no use of the prior measure.  Note that the optimal prediction
runs yield good results when the standard deviation of the values of
the collective variables is as large as 50 percents of their mean (the
standard deviation of the pointwise values of the solutions is much
larger still).

\section{Pseudo-spectral optimal prediction for a model nonlinear problem}
\label{sec:dealiasing}
\setcounter{equation}{0}
\setcounter{figure}{0}
\setcounter{table}{0}

In the present section we consider a discrete version of the same
Schr\"odinger equation as above; our goal is to show explicitly how the
information in the prior measure improves the accuracy of an
underresolved nonlinear calculation. One of the striking facts shown
by the example is that first-order optimal prediction is useful in
nonlinear problems.  In view of the rapid convergence of the
conditional expectations of discrete problems to their continuum
limits, as displayed in the previous section, we shall be content with
the discrete problem.  Specifically, we shall contrast the solution of
a discrete problem by a pseudo-spectral method that takes no
cognizance of the prior measure with a closely related optimal
prediction scheme with Fourier kernels, and show how the information
in the prior measure improves the predictions.  In the present
section, a complex-function formalism turns out to be more
transparent, and we therefore slightly change the notations.  We
consider a set of complex ordinary differential equations which is a
formal discretization of our Schr\"odinger equation,
(\ref{eq:complex.schrodinger}),
\begin{equation}
\imath \pd{u_j}{t} = -\frac{u_{j-1} - 2 u_j + u_{j+1}}{h^2}
+ \frac{1}{4} \Brk{3 \Abs{u_j}^2 u_j + (u_j^*)^3},
\label{nse:dynamics}
\end{equation}
with $j=1,\ldots,n$ and $h=2\pi/n$. Equations (\ref{nse:dynamics}) are
the Hamilton equations of motion derived from the following
Hamiltonian,
\begin{equation}
{\cal H}(u) = \frac{1}{2} \sum_{j=1}^n
\BRK{ \frac{\Abs{u_{j+1}-u_j}^2}{h} + \frac{h}{16}
\Brk{u_j^4 + 6 \Abs{u_j}^4 + (u_j^*)^4} } 
\label{nse:H}
\end{equation}
(see section $4$ above). One can readily verify that this is the same
discrete Hamiltonian as before. The prior measure is the canonical
measure whose density is
\begin{equation}
f(u) = \frac{1}{Z} e^{-{\cal H}(u)}.
\label{nse:prior}
\end{equation}
Note that we write $u$ without boldface, as we did in the case of
functions, but not as we did for vectors; this is done by analogy with
the previous sections on the Schr\"odinger equation.  (Those who look
at our earlier papers \cite{ckk1},\cite{ckk2},\cite{ckk3} will notice
that the measure here differs from the measure used there by a factor
$h$ in the exponent, and also that we used a different spatial
period.)  Equation (\ref{nse:prior}) requires an explanation when the
variables $u_j$ are complex.  A set of $n$ complex variables can be
treated as a set of $2n$ real variables.  However, when the real and
imaginary parts $q_i,p_i$ of each variable $u_i$ are independent and
their means are zero, (as they are here), one can write:

\bea
\lefteqn{P(s_1 < q_1 \le s_1 + ds_1, \dots, r_1 < p_1 \le r_1+dr_1,\dots, 
r_n < p_n \le r_n+dr_n) =} \nonumber \\
&& = F({\mathbf z}) \, dz_1\dots dz_n
= Z^{-1} \exp \brk{- \half (\bfz,A\bfz^*)} dz_1\dots dz_N,
\label{complexA}
\eea
where $z_i=s_i+\imath r_i, dz_i=ds_i dr_i$, and the matrix $A$ is
hermitian.  When ${\cal H}$ is a quadratic function of the vector $u$,
equation (\ref{nse:prior}) defines $A$ in the complex case.  

The prior measure (\ref{nse:prior}) being non-Gaussian, we proceed as
above and partition the Hamiltonian into a quadratic part plus a
perturbation, ${\cal H} = {\cal H}^0 + {\cal H}^1$. To make the
example amenable to analysis we keep only the leading term in this expansion; we
already pointed out that the leading term contains a contribution of
the nonlinear terms in the equation and that the partition takes into account higher order terms in the expansion.  As explained above, there is here
only one relevant partitioning parameter, $m_0$.  Thus we approximate
the probability density by $f(u) = Z_0^{-1} e^{-{\cal H}^0(u)}$,
where
\begin{equation}
{\cal H}^0(u) = \frac{1}{2} \sum_{j=1}^n
\BRK{ \frac{\Abs{u_{j+1}-u_j}^2}{h} + h m_0^2  \Abs{u_i}^2}.
\label{nse:H0}
\end{equation}
We work here with $n=32$.  In this finite dimensional case, we
rederived the optimal value of $m_0$ as follows: We calculated the
two-point correlation function $\Average{u_i u_j^*}^0$ for the measure
(\ref{nse:prior}) by a Monte-Carlo method, and then found the value of
$m_0$ that best reproduces this correlation function by minimizing the
mean-square difference between this correlation function and the
correlation function produced by (\ref{nse:H0}). This yielded
$m_0=1.055$, confirming the value obtained above.  The procedure used
here minimized the difference between the full measure and its
quadratic piece while the more general procedure of the previous
sections minimizes only the first few terms in an expansion; it is
comforting that the results agree. For readers of our earlier papers
\cite{ckk1,ckk2,ckk3}, we point out that the present procedure differs
from the ``Gaussianization'' proposed there by using an analytical
expression for the approximate measure, whereas in the previous
publications we needed to store the full covariance matrix.
Furthermore, the present construction produces a first term in a
systematic expansion.

\begin{figure}
 \centerline{\scalebox{0.5}{\includegraphics{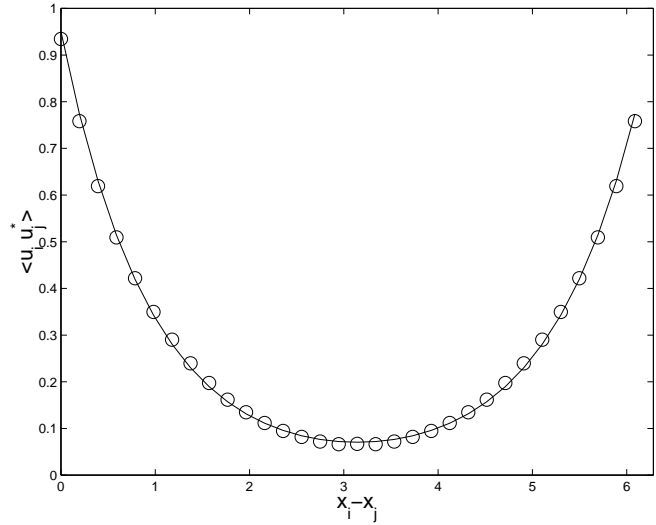}}}
\caption{{\sf Comparison between the two-point correlation function
$\Average{u_i \,u^*_j}$, as computed by a Monte Carlo sampling
procedure, and its approximation (\ref{nse:Cij}) based on the Gaussian
measure (\ref{nse:H0}), with $m_0=1.05$. The number of points is
$n=32$. }}
\label{fig:gaussianization}
\end{figure}

In Figure \ref{fig:gaussianization} we compare the exact two-point
correlation function $\Average{u_i\,u^*_j}$ obtained by a Monte Carlo
sampling (open circles) to the Gaussian approximation (solid line),
\begin{equation}
C_{ij} = \Average{u_i\,u^*_j}^0 = 
\frac{1}{\pi} \sum_{k=1}^n \frac{e^{\imath k(x_i-x_j)}}%
{\frac{4}{h^2} \sin^2 \frac{h k}{2} + m_0^2},
\label{nse:Cij}
\end{equation}
derived from (\ref{nse:H0}). The discrepancy is negligible compared
with the statistical uncertainty in the sampling procedure.

The formulas for calculating conditional expectations 
(\ref{constrainedMean})--(\ref{constrainedCovariance}) can be
generalized to deal with complex functions. Assume we have $N$ conditions 
of the form,
\[
g_{\alpha i} u_i = V_\alpha,
\]
where the $g_{\alpha i}$ are the complex entries of the $N\times n$
kernel matrix $G$ and the entries of the vector  $V_\alpha$ are the values of the
$N$ collective variables defined
by $G$. The formulas for the conditional means and variances of the
complex vector $u_i$ generalize equations 
(\ref{constrainedMean}) and (\ref{constrainedCovariance}); 
the conditional average of $u_i$ is 
\[
\Average{u_i}^0_V = q_{i \alpha} V_\alpha;
\]
the $q_{i \alpha}$ are the entries of the matrix $Q = (C G^\dagger)(G C G^\dagger)^{-1}$,
where $C$, whose elements are defined in (\ref{nse:Cij}), is the
matrix that approximates the inverse of the matrix $A$ defined in
(\ref{complexA}). The dagger denotes an adjoint (hermitian transposed) matrix.  The conditional covariance matrix has entries
\[
\Cov{u_i}{u^*_j}^0_V = \Brk{C - (C G^\dagger)(G C G^\dagger)^{-1}(G C)}_{ij}.
\]

We now make a special choice of kernels $G$: We pick them so that
what is known at time $t=0$ is a set of Fourier modes,
fewer than are necessary to specify the solution completely.
This makes the $g_{\alpha i}$ complex exponentials:  
\[
g_{\alpha i} = \frac{1}{n} \exp\brk{-\imath K_\alpha x_i}.
\]
If the number of conditions $N$ is even, the 
$K_\alpha$ take the values $-\frac{N}{2}+1,\ldots,\frac{N}{2}$.  
Note the following property of the resulting matrix $G$: 
\begin{equation}
GG^\dagger=\frac {1}{n}I,
\label{identity}
\end{equation} 
where $I$ is the $N\times N$ identity.  Spectral
variables are particularly convenient here because, if $u$ is expanded
in Fourier series, and this series is substituted into the formula for
the Hamiltonian ${\cal H}$, the result is a sum of squares of Fourier
coefficients; the prior measure is ``diagonal in Fourier space.'' This
can also be deduced from the fact that the measure is Gaussian and
translation-invariant in the variable $x$. In particular, after $u$ is
expanded in Fourier series, the matrices $A, C$ are both diagonal.  A
short calculation yields:
\[
Q = n\,G^\dagger,
\]
from which follows
\begin{equation}
\Average{u}^0_V = n\,G^\dagger V,
\label{nse:mean}
\end{equation}
and 
\begin{equation}
\Cov{u_i}{u^*_j}^0_V = \Brk{C - n\,G^\dagger G C}_{ij}.
\label{nse:cov}
\end{equation}
Equation (\ref{nse:mean}) is the interpolation formula for the first
moment of $u$. Recall that $V=G u$, hence $\Average{u}_V=n\,G^\dagger
G u$. The operator $n\,G^\dagger G$ is a Galerkin projection operator
which projects any vector onto the vector space spanned by the range
of $G$. Thus, the mean of $u$ produced by our regression formula
equals the mean obtained from a simple Fourier series that uses the
known coefficients $V_{\alpha}$ (we shall call this Fourier series the
``Galerkin reconstruction''). This is as it should be: In a
translation-invariant Gaussian measure the Fourier components are
mutually independent; the knowledge of the first few does not
condition the next ones, whose expected value is therefore zero.
However, the measure does contain information about the higher moments
of the higher Fourier coefficients, and this is important in a
nonlinear problem.

\begin{figure}
 \centerline{\scalebox{0.5}{\includegraphics{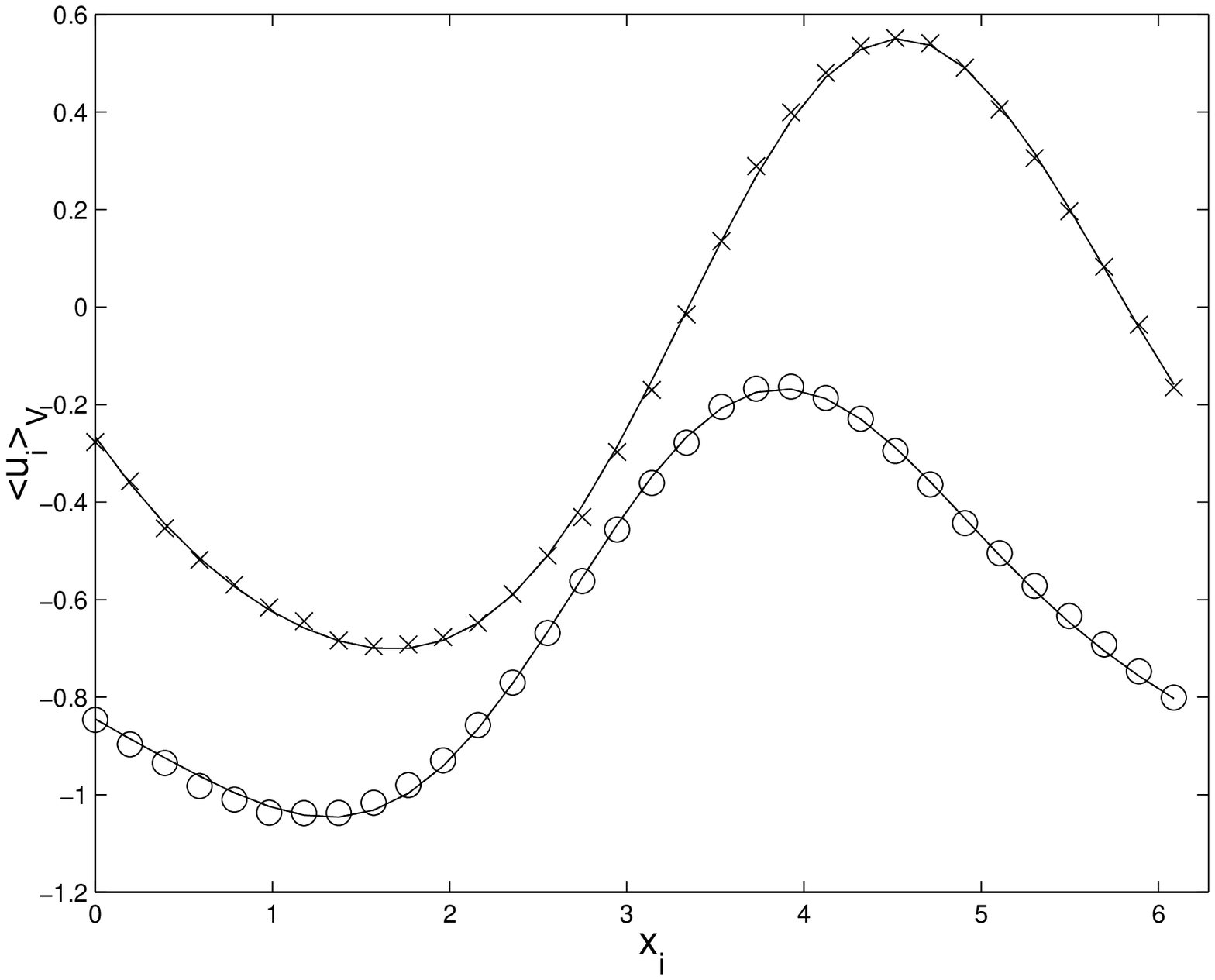}}}
 \centerline{\scalebox{0.5}{\includegraphics{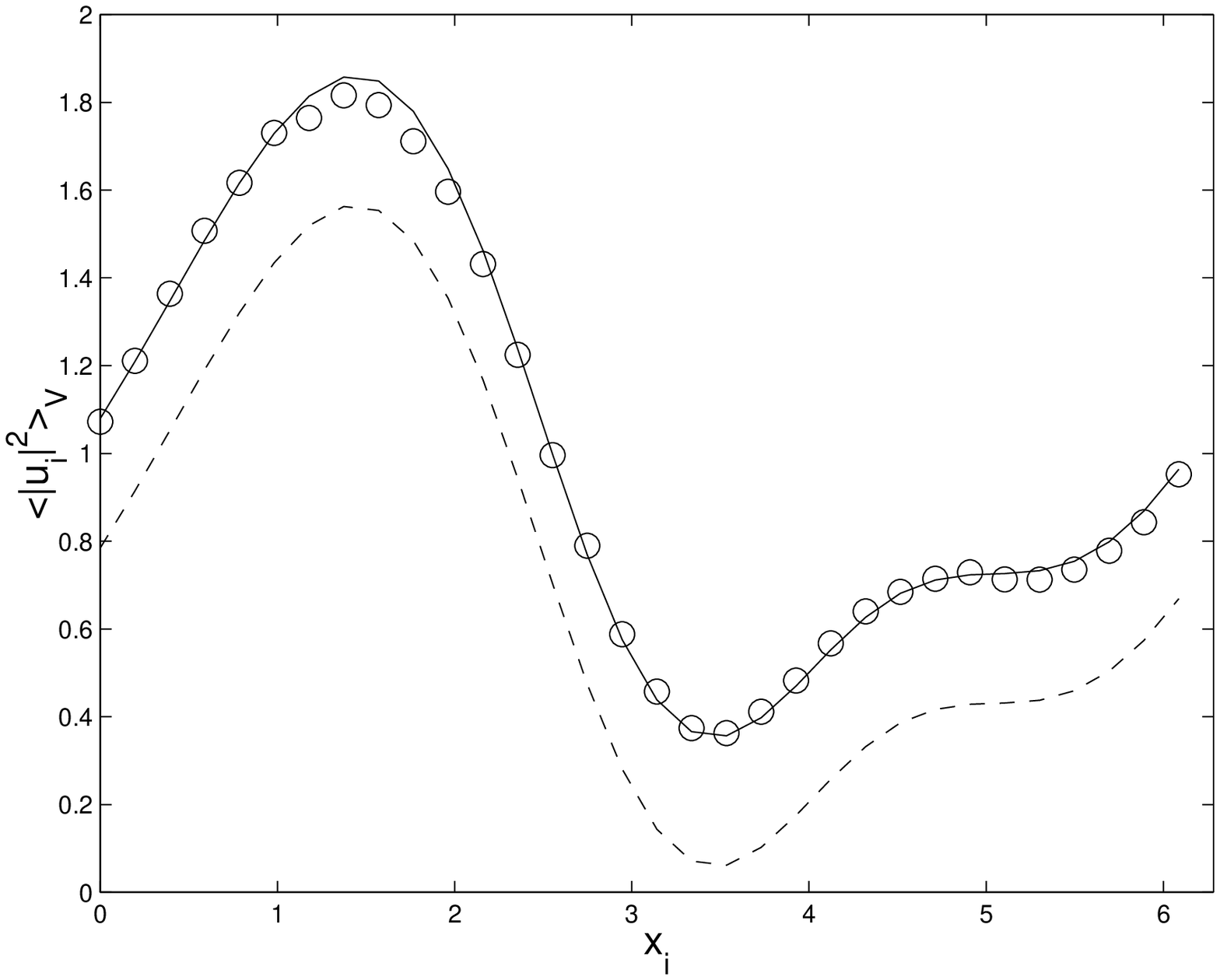}}}
 \centerline{\scalebox{0.5}{\includegraphics{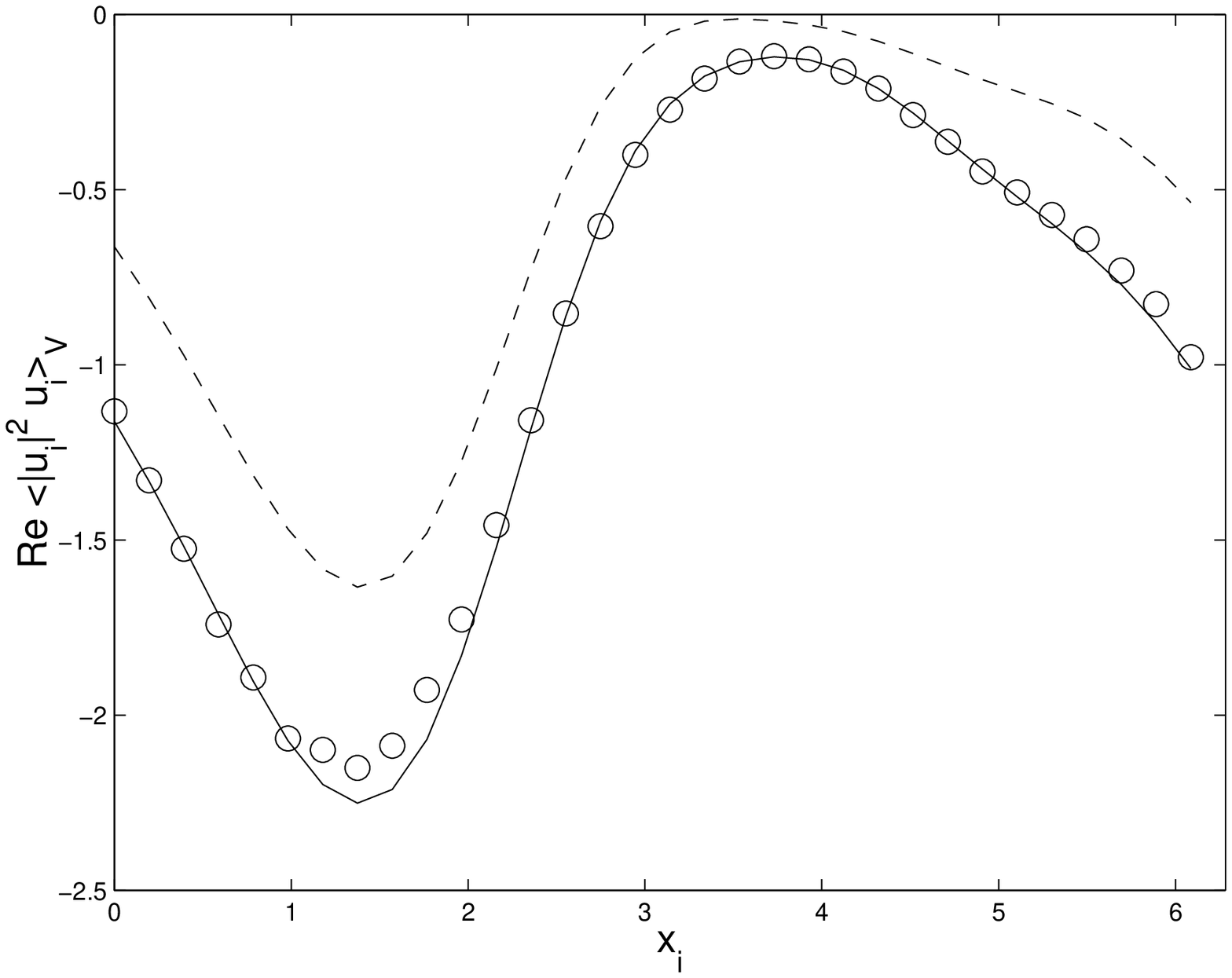}}}
\caption{{\sf Conditional moments of $u$: comparison between true 
values (symbols), a Galerkin reconstruction (dashed
lines), and regression (solid lines). (a) The conditional
expectation $\Average{u_i}_V$; the circles represent the real part and
the crosses represent the imaginary part.  (b) The variance
$\Average{u_i^2}_V$. (c) The third moment
$\mbox{Re}\,\Average{\Abs{u_i}^2 u_i}_V$.  The number of points is
$n=32$, and the number of Fourier modes that are assumed known is $N=4$.}}
\label{fig:interpolation}
\end{figure}

The difference between the Galerkin reconstruction and the regression
used in the optimal prediction of moments is demonstrated in Figure
\ref{fig:interpolation} . The three graphs depict the first three
moments, calculated (i) by a Monte Carlo sampling of the measure
(\ref{nse:prior}) conditioned by $N=4$ known collective variables
(symbols), a Galerkin reconstruction (dashed lines), and our
regression formulas (solid lines). We chose a system size of $n=32$
with $N=4$ resolved (i.e., known) Fourier modes. For the first moment,
$\Average{u_i}_V$ (calculated by Monte-Carlo), the Galerkin
reconstruction and the regression are close to each other. For the
second and third moments the regression, which is the core of optimal
prediction, is close to the truth, as revealed by the Monte-Carlo
runs; the error is smaller than the statistical uncertainty in the
sampling. The Galerkin reconstruction, on the other hand, deviates
significantly from the truth. These graphs demonstrate the importance
of the prior measure, which contains information about the mean
squares of the unresolved Fourier components.

We next derive the optimal prediction scheme for our model problem
with Fourier kernels.  Given a set of Fourier modes, $V_\alpha$, we
replace the right-hand side of (\ref{nse:dynamics}) by its conditional
average, and multiply the result by the kernel matrix $G$ . Thus,
\[
\imath \deriv{V_\alpha}{t} = g_{\alpha i} 
\Average{- \frac{u_{i+1}-2 u_i+u_{i-1}}{h^2} + 
\frac{1}{4} \Brk{3 \Abs{u_i}^2 u_i + (u^*_i)^3}}^0_V.
\]
Substituting the regression formulas (\ref{nse:mean}) and
(\ref{nse:cov}) and using Wick's theorem and (\ref{identity}), we
obtain the following set of $N$ equations:
\begin{equation}
\imath \deriv{V_\alpha}{t} =
\frac{4}{h^2} \sin^2 \frac{K_\alpha h}{2} V_\alpha + 
\frac{1}{4} \sum_{\beta,\gamma,\epsilon} 
\Brk{3 V_\beta V^*_\gamma V_\epsilon
\,\delta_{\alpha,\beta-\gamma+\epsilon} +
V^*_\beta V^*_\gamma V^*_\epsilon
\,\delta_{\alpha,-\beta-\gamma-\epsilon}} + \frac{6}{4}c\,V_\alpha,
\label{nse:effective}
\end{equation}
where
\[
c = \Brk{(I - n\,G^\dagger G)C}_{ii}
\qquad\mbox{(no summation on $i$)}
\]
is a constant (the right-hand side is independent of $i$). Note that
the last term comes from the evaluation of the nonlinear terms by
Wick's theorem.  The structure of equation (\ref{nse:effective}) is
enlightening: The first two terms on the right-hand side are precisely
the Galerkin approximation for the evolution of a subset of Fourier
modes; they constitute a pseudo-spectral approximation of the
equations of motion. The third term, which is linear in $V$,
represents information gleaned from the prior measure.  The nice
feature of this example is the sharp separation between the
contribution from the resolved degrees of freedom and the contribution
of the ``subgrid'' degrees of freedom, which happens to simply
``renormalize'' the linear part of the evolution operator.

\begin{figure}
 \centerline{\scalebox{0.5}{\includegraphics{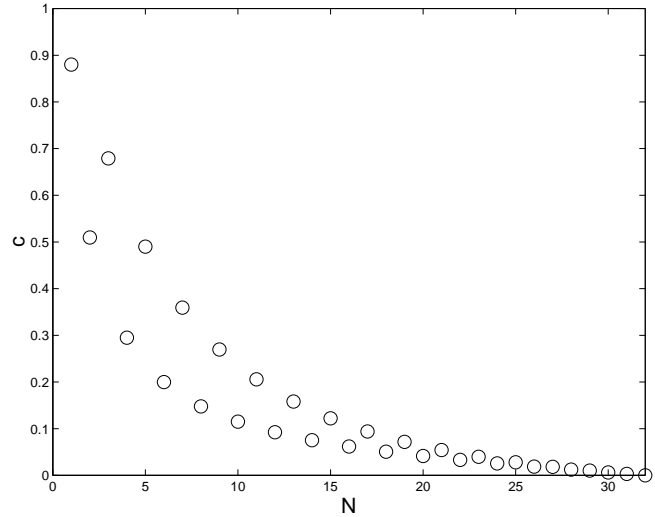}}}
\caption{{\sf The parameter $c$ as a function of the number $N$ of
resolved Fourier modes, for $n=32$.}}
\label{fig:CofN}
\end{figure}

One is of course interested in knowing how large is the extra term
that makes up the entire difference between the optimal prediction
scheme and a standard pseudo-spectral scheme; this difference is
proportional to the coefficient $c$. In Figure \ref{fig:CofN} we plot
the value of $c$ as function of $N$ for $n=32$. As expected, $c$ is
larger when the number of resolved degrees of freedom is smaller , and
vanishes when $N=n$, i.e., when the system is fully resolved. The
oscillations in this graph result from the alternation between odd and
even numbers of Fourier modes. Figure \ref{fig:CofN} demonstrates that
optimal prediction is consistent: As the number of collective
variables increases, its predictions converge to those of a resolved
calculations (as indeed should be obvious from the derivation); when
the number of collective variables is small, the corrections due to
optimal prediction are substantial.

\begin{figure}
\centerline{\scalebox{0.5}{\includegraphics{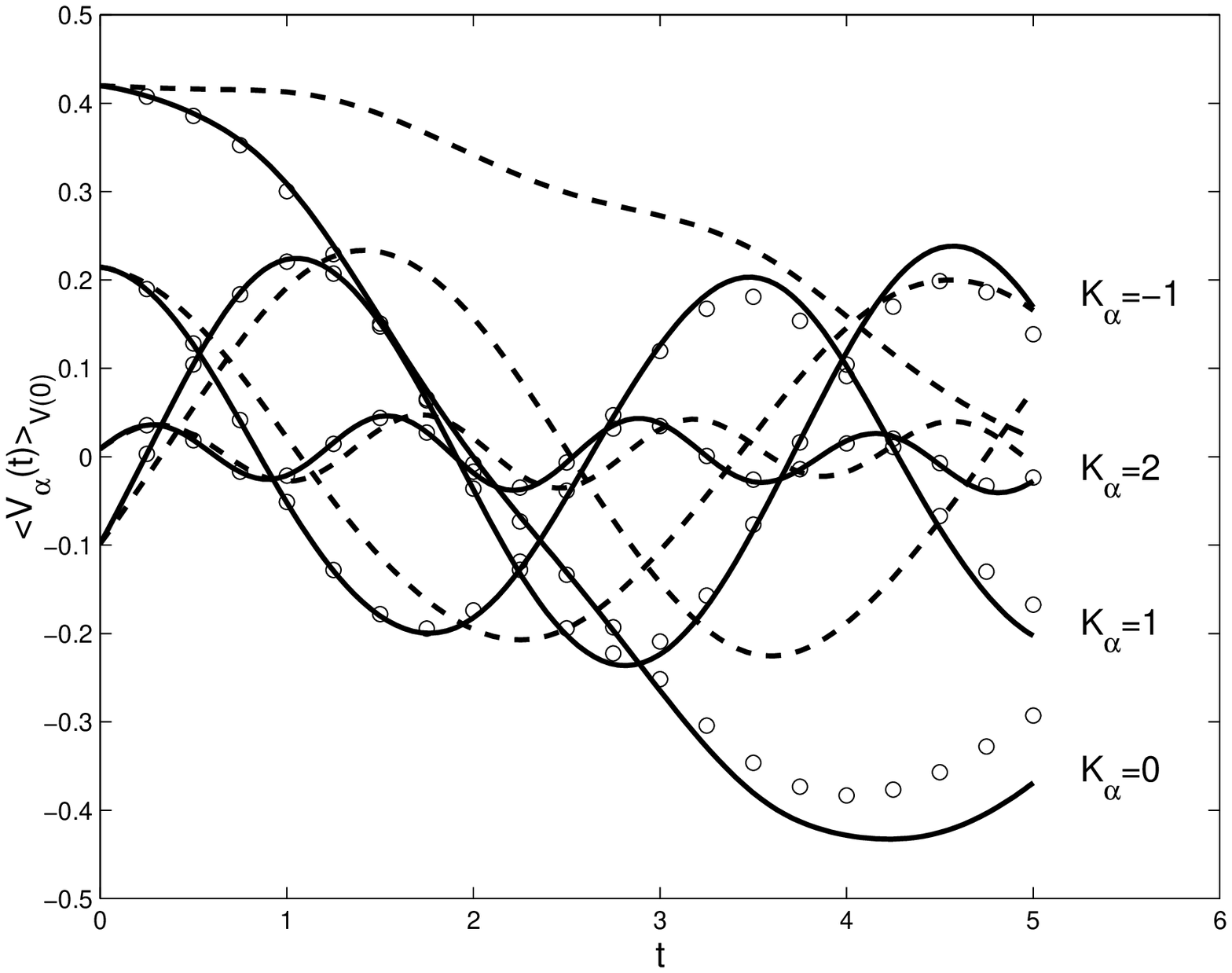}}
            \scalebox{0.5}{\includegraphics{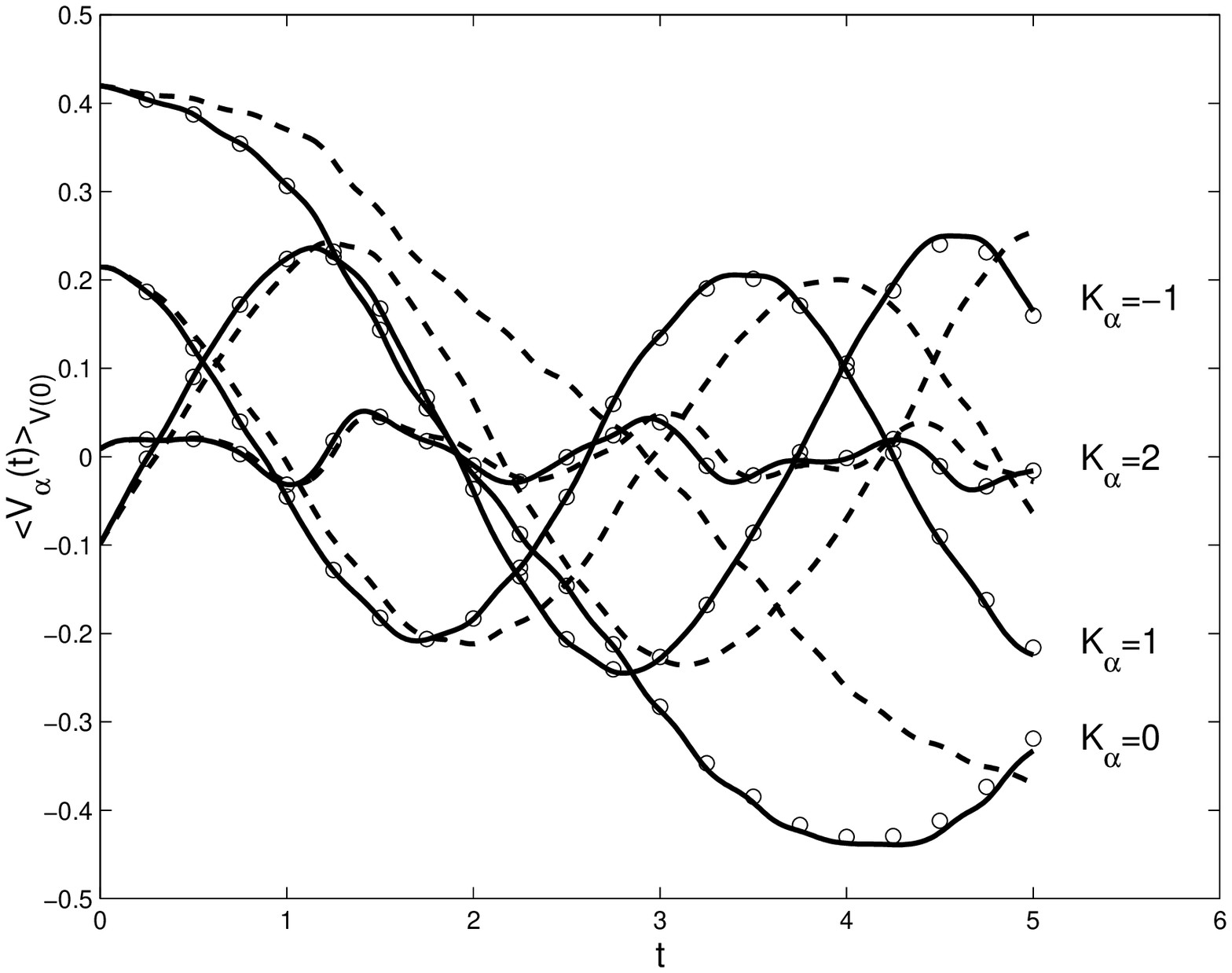}}}
\caption{{\sf Time evolution of the real part of the $4$ lowest
Fourier modes for $n=32$. The circles represent average values over an
ensemble of $10^4$ states; the dashed lines result from the Galerkin
approximation; the solid lines result from our optimal prediction
scheme (\ref{nse:effective}).  Results are presented for (a) $N=4$ and
(b) $N=8$.}}
\label{fig:time}
\end{figure}

In Figure \ref{fig:time} we compare the time evolution of the first 4
Fourier coefficients predicted by our optimal prediction scheme
(\ref{nse:effective}) (solid lines), the time evolution of these modes
predicted by a Galerkin scheme (dashed lines), and their exact mean
evolution obtained by sampling $10^4$ states from the conditional
measure, evolving them in time, and averaging the first Fourier
components (circles).  We do this with $N=4$ and $N=8$; with $N=8$ we
display only $4$ modes even though $8$ are calculated, in the interest
of clarity. To make the calculations with different values of $N$
comparable, we pick the initial values of the collective variables in
the way suggested by the analysis of Hald \cite{hald:optimal}: We
sample an initial function $u$ from the invariant measure, and then
calculate values of the collective variables by performing the
summations $g_{\alpha i}u_i$.  The graph shows that in each case the
simple Galerkin calculation deviates immediately from the true
solution, while the optimal prediction remains accurate. The
calculations show that optimal prediction improves the accuracy
compared to a Galerkin calculation; the time during which the optimal
prediction remains accurate increases with increasing $N$; we know
from (\ref{nse:effective}) that the cost of optimal prediction in this
problem is small.

\section{Conclusions}

We have exhibited the value of the statistical information used in
optimal prediction for the solution of an underresolved nonlinear
problem, and we have shown that perturbation theory provides a
ready-made machinery for applying the ideas of optimal prediction to
problems where the invariant measure is non-Gaussian.

The first-order implementation of optimal prediction with a fixed, small, number
of conditions breaks down after a finite time; the time for which it is valid
increases as the temperature increases;
the temperature determines the variance of the invariant measure and thus the uncertainty
in the system.
There are two ways to improve the prediction: Go to more
sophisticated prediction theory, as outlined in section 2, or increase
the number of collective variables. We have demonstrated the power of the
second alternative; the first alternative will be
explored in later publications.

Many aspects of the algorithms presented here require further work.
The closure by means of a fixed number of affine conditions is only a
first step; other closure schemes will be investigated. In particular,
optimal prediction fits within the framework of irreversible
statistical mechanics, whose apparatus can be brought into use for
finding closure schemes.  More powerful versions of perturbation
theory can be readily used.  A careful perusal of our final example
shows that there is a great advantage in using orthogonal kernels; in
the interest of pedagogy, we have not used this possibility in the
discussion of perturbation theory, and we have yet to explore
orthogonal bases other than Fourier bases.

An inspection of the formulas derived by perturbation theory shows
that though we assumed a knowledge of an invariant measure, all that
is finally used is a set of moments; this is the opening for applying
optimal prediction methods in problems where less than a full
invariant measure is known.

We assume that all these issues will be handled as we progress to more
complicated problems, in more dimensions, with dissipation (requiring
a careful modeling of effective Hamiltonians), and with general
boundary conditions. We shall explore problems with these additional
features in future publications.


\subsection*{Acknowledgments} 
The authors would like to thank Prof. G.I. Barenblatt,
Mr. T. Burin des Roziers, Dr. A. Gottlieb, Prof. O. Hald, Dr. A. Kast, 
Prof. D. Kessler, Dr. I. Kliakhandler,  and Prof. P. Marcus 
for helpful discussions and comments.  This work was supported in part
by the Applied Mathematical Sciences subprogram of the Office of
Energy Research, U.S. Department of Energy, under contract
DE--AC03--76--SF00098, and in part by the National Science Foundation
under grant DMS94--14631.


\begin{thebibliography}{99}

\bibitem{bennett} Bennett, A.F.,
{\em Inverse methods in physical oceanography}, Cambridge Monographs
on Mechanics and Applied Mathematics, Cambridge, 1992

\bibitem{chorin} Chorin, A.J., 
{\em Accurate evaluation of Wiener integrals},  Math. Comp. 27 (1973), 1-15.

\bibitem{chorin2} Chorin, A.J., 
{\em Vortex methods}, Les Houches Summer School of Theoretical Physics,
59, Elsevier (1995), 67-108.

\bibitem{ckk1} Chorin A.J., Kast A.P., Kupferman R., 
{\em Optimal prediction of underresolved dynamics}, Proc. Nat. Acad. Sci. USA,
{\bf 95}, (1998), 4094-4098.

\bibitem{ckk2} Chorin A.J., Kast A.P., Kupferman R.,
{\em Unresolved computation and optimal prediction},
Comm. Pure Appl. Math., in press, 1999.

\bibitem{ckk3} Chorin A.J., Kast A.P., Kupferman R.,
{\em On the prediction of large-scale dynamics using unresolved computations},
Contemp. Math., in press, 1999.

\bibitem{evans}
Evans, D., Morriss, G.,
{\em Statistical Mechanics of Non-equilbrium Liquids},
Academic Press, NY, 1990.

\bibitem{fetter-walecka:quantum} Fetter A.L., Walecka J.D.,
{\em Quantum theory of many particle systems},
McGraw-Hill, 1971.

\bibitem{freedman} Freedman, F., 
{\em Brownian motion and diffusion}, Springer, NY, 1983

\bibitem{goldstein} Goldstein, H., 
{\em Classical mechanics}, Addison Wesley, Reading, Mass, 1980.

\bibitem{gottlieb} Gottlieb, A.D., 
{\em First order optimal prediction}, in preparation.

\bibitem{hald:optimal} Hald O.H.,
{\em Optimal prediction and the Klein-Gordon equation},
Proc. Nat. Acad. Sci. USA 96 (1999), 4774-4779.

\bibitem{kleinert:gauge} Kleinert H.,
{\em Gauge fields in condensed matter}, Vol. I,
World-Scientific, Singapore, 1989.

\bibitem{ma} Ma S.-K.,
{\em Modern theory of critical phenomena},
W.A. Benjamin, Massachusetts, 1976.

\bibitem{mckean} MacKean, H.,  
{\em Statistical mechanics of nonlinear wave equations IV: Cubic
Schr\"odinger equations}, Comm. Math. Phys. 168, (1995), 479-491.

\bibitem{mori}
Mori, H., 
{\em Transport, collective motion, and Brownian motion}
Prog. Th. Phys., 33, (1965), 423-455.

\bibitem{morokoff} Morokoff, W., 
{\em Generating quasi-random paths for stochastic processes}, 
SIAM Rev., 40, (1998), 765-788.

\bibitem{mandl} Mandl, F., 
{\em Introduction to quantum field theory}, Interscience, NY, 1959.

\bibitem{lax} Mock M., Lax P., 
{\em The computation of discontinuous solutions of linear
hyperbolic equations}, Comm. Pure Appl. Math., 31, (1978), 423-430.

\bibitem{onsager} Onsager, L., 
{\em Statistical hydrodynamics}, Nuovo Cimento, Suppl. to Vol. 6, 
(1949), 279-287.

\bibitem{papoulis} Papoulis, A., 
{\em Probability, random variables and stochastic processes},
McGraw-Hill, New York, 1991.

\bibitem{proudman} Proudman, I., Reid, W.,  
{\em On the decay of normally distributed and homogeneous turbulent
velocity fields}, Phil. Trans. Roy. Soc. A247, (1954), 163-189.

\bibitem{ramond} Ramond, P., 
{\em Field theory, a modern primer}, Benjamin, Reading, Mass, 1981

\bibitem{scotti}
Scotti, A., Meneveau, C., {\em Fractal model for coarse-grained
partial differential equations}, Phys. Rev. Lett. 78, (1997), 867-870.

\bibitem{vaillant}
Vaillant, O., PhD thesis, INRIA Sophia-Antipolis, 1999. 

\bibitem{zwanzig}
Zwanzig, R., 
{\em Problems in nonlinar transport theory},
Systems Far From Equilibrium, 
L. Garrido (ed),
Interscience, New York, 1961, 198-221.

\end{thebibliography}
\end{document}